\documentclass[a4,12pt]{article}

\usepackage{amssymb,amsmath}
\usepackage{xypic}
\usepackage{graphicx}

\newtheorem{defn}{Definition}[section]
\newtheorem{proposition}[defn]{Proposition}
\newtheorem{corollary}[defn]{Corollary}
\newtheorem{rem}[defn]{Remark}
\newtheorem{exm}[defn]{Example}
\newtheorem{lemma}[defn]{Lemma}
\newtheorem{theorem}[defn]{Theorem}
\newtheorem{notat}[defn]{Notation}
\newtheorem{newpar}[defn]{}

\newtheorem{xdefn}{Definition.}
\newtheorem{xproposition}{Proposition.}
\newtheorem{xcorollary}{Corollary.}
\newtheorem{xrem}{Remark.}
\newtheorem{xexm}{Example.}
\newtheorem{xlemma}{Lemma.}
\newtheorem{xtheorem}{Theorem.}
\newtheorem{xnotat}{Notation.}
\newtheorem{xnewpar}{\it}
\newtheorem{xproof}{{\it Proof. }}
\newtheorem{xproofof}{{\it Proof}}

\newenvironment{definition}{\begin{defn}\em}{\end{defn}}

\newenvironment{example}{\begin{exm}\em}{\end{exm}}

\newenvironment{proof}{\begin{xproof}\em}{\end{xproof}}

\newenvironment{newparagraph*}[1]{\begin{xnewpar}\hspace*{-1.5mm}{#1}. \rm}{\end{xnewpar}}

\newenvironment{definition*}{\begin{xdefn}\em}{\end{xdefn}}
\newenvironment{remark*}{\begin{xrem}\em}{\end{xrem}}
\newenvironment{example*}{\begin{xexm}\em}{\end{xexm}}
\newenvironment{notation*}{\begin{xnotat}\em}{\end{xnotat}}
\newenvironment{proposition*}{\begin{xproposition}}{\end{xproposition}}
\newenvironment{corollary*}{\begin{xcorollary}}{\end{xcorollary}}
\newenvironment{lemma*}{\begin{xlemma}}{\end{xlemma}}
\newenvironment{theorem*}{\begin{xtheorem}}{\end{xtheorem}}

\newcommand\N{\mathbb{N}}
\def\qed{\hspace{0.3cm}{\rule{1ex}{2ex}}}
\newcommand\V{\bigvee}
\newcommand\Q{\mathrm{Max\,}}
\newcommand\ie{i.e.}
\newcommand\eg{e.g.}
\newcommand\then{\&}
\newcommand\seq{{\rm seq}}
\newcommand\pwset[1]{{\cal P}({#1})}
\newcommand\st{\mid}
\newcommand\goto[1]{\langle{#1}\rangle}
\newcommand{\bprop}[2]{{#2}^{(#1)}}
\newcommand{\wprop}[2]{{\langle \bprop{#1}{#2}\vert}}
\newcommand{\rprop}[2]{{\vert \bprop{#1}{#2}\rangle}}
\newcommand{\prop}[2]{\bprop{#1}{#2}}
\newcommand{\entails}{\;\vdash\;}
\newcommand{\true}{\mathit{true}}
\newcommand{\false}{\mathit{false}}
\newcommand\cantor{C}
\newcommand\ann{\mathrm{ann}}

\newcommand\End{\mathcal{Q}}
\newcommand\Rel{\mathcal{Q}}
\newcommand\cf{\textrm{cf.}}
\newcommand\pent{\mathbf{Pen}}
\newcommand\penq{\mathrm{Pen}}
\newcommand\pens{\pent_s}
\newcommand\pena{A}
\newcommand\Seq{\seq^\sharp} 

\begin{document}

\title{A noncommutative theory of Penrose tilings\thanks{Research supported in part by FEDER
and FCT/POCTI/POSI through the research units funding program and grant POCTI/1999/MAT/33018, and by CRUP and British Council through Treaty of Windsor grant No.\ B-22/04.}}
\author{Christopher J.~Mulvey$^1$ and Pedro Resende$^2$
\vspace*{2mm}\\ \small\it $^1$School of Mathematical Sciences, University
of Sussex,\vspace*{-2mm}\\ \small\it Falmer, Brighton, BN1 9QH, United
Kingdom\\ \small\it $^2$Departamento de Matem{\'a}tica, Instituto
Superior T{\'e}cnico, \vspace*{-2mm}\\ \small\it Av. Rovisco Pais,
1049-001 Lisboa, Portugal}

\date{~}

\maketitle

\vspace*{-1,3cm}
\begin{abstract}
Considering quantales as generalised noncommutative
spaces, we address as an example a quantale $\penq$ based on the Penrose
tilings of the 
plane. We study in general the representations of involutive 
quantales on those of binary relations, and show that in the case of $\penq$
the algebraically irreducible representations provide a 
complete classification of the set of Penrose tilings from which its
representation as a quotient of Cantor space is recovered.
\vspace{0.1cm}\\ \textit{Keywords:} quantale, noncommutative space, noncommutative logic, Penrose tiling, relational representation, relational point.
\vspace{0.1cm}\\ 2000 \textit{Mathematics Subject
Classification}: Primary 06F07, 52C23; Secondary 03B60, 05B45, 18B99, 
46L85, 54A05, 58B34, 81P10, 81R60, 81S99.
\end{abstract}

\maketitle

\section{Introduction}\label{sec:introduction}

The concept of \emph{quantale}, introduced in~\cite{Mu86} to investigate the spectrum of noncommutative C*-algebras, generalises that of locale, or ``point free space" (\cf\ \cite{Johnstone}), in a way that brings together the noncommutative topology of~\cite{GiKu71} with ideas stemming from the constructive formulation of classical results in the theory of commutative C*-algebras~\cite{BaMu00a,BaMu00b,BaMu:prep}. In particular, there is a functor from unital C*-algebras to the category of unital involutive quantales which is a complete invariant: to each C*-algebra $A$ it assigns the quantale $\Q A$ of all the closed linear subspaces of $A$ (its operator spaces~\cite{Pisier}). This functor was proposed in~\cite{Mu89} as a noncommutative generalisation of the classical maximal spectrum, and has been subsequently studied in various papers~\cite{Kr02:phd,KrPeReRo,KrRe,Mu02,MuPe01,MuPe02,PaRo00,Re:2-forms}.

Another motivation for the notion of quantale in~\cite{Mu86} was that presentations of quantales by generators and relations can be regarded equivalently as theories in a noncommutative logic whose conjunction $\then$ is related to the arrow of time: $a\then b$ is to be read as ``$a$ and \emph{then} $b$", emphasising the fact that $a\then b$ is in general different from $b\then a$, and $a\then a$ is different from $a$. This logical view is analogous to the Lindenbaum construction that allows theories in propositional geometric logic to be considered as presentations of locales by generators and relations, which in particular was exploited when proving in an arbitrary Grothendieck topos that the category of commutative unital C*-algebras is equivalent to the category of compact completely regular locales~\cite{BaMu00a,BaMu00b,BaMu:prep}, thus generalising the classical Gelfand duality, which (in the presence of the Axiom of Choice) establishes an equivalence between (the dual of) the category of commutative unital C*-algebras and the category of compact Hausdorff spaces.

This logical view also suggests the possibility of describing quantum systems using quantales determined by noncommutative theories~\cite{Mu86}, although this idea has not so far been pursued in the context of quantum physics: quantales have been mentioned in physics, but not in this way (see, for instance, \cite{AmCoSt98,CoSt99}). However, the logical view was one of the influences behind the introduction of the quantale $\Q A$, in the sense that a presentation by generators and relations may be obtained as a generalisation of that of the locale $\Q A$ considered in the case of a commutative C*-algebra $A$~\cite{BaMu00a,BaMu00b,BaMu:prep}. It has also been used extensively in computer science when studying concurrent systems~\cite{AbVi93,Re99,Re00,Re01,Re02:jpaa,ReVi}.

Since quantales may also be considered as semigroups in the closed category of sup-lattices~\cite{JoTi84}, a natural way of studying them is to look at their actions by endomorphisms on sup-lattices, that is, their right or left modules, or, equivalently, their representations. In particular, the notion of \emph{algebraically irreducible representation} of an involutive quantale that was introduced in~\cite{MuPe01} provides a definition of ``point" of a quantale that coincides, in the case of locales, with the usual notion of point. 
More precisely, we take a \emph{point} to be (the equivalence class of) any such representation. 

In the case of the quantale $\Q A$, the points provide a complete classification of the irreducible representations of a unital C*-algebra $A$, in the sense of a bijective correspondence between the points of $\Q A$ and the equivalence classes of irreducible representations of $A$. In fact, the result actually proved in~\cite{MuPe01} is that this is the case exactly for those points for which there exists a pure state of the C*-algebra of which the kernel is mapped properly by the representation, it being conjectured that this is indeed the case for every algebraically irreducible representation of $\Q A$. The sup-lattices on which $\Q A$ acts are in this case shown to be those of projections of the Hilbert spaces on which $A$ is irreducibly represented. In contrast, representations on complete atomic Boolean algebras are called \emph{relational}. They are important in computer science, where they can be identified with concurrent systems~\cite{Re00,Re01,Re02:jpaa}.

In order to obtain better insight into the relation between quantales and C*-algebras, in this paper we look, albeit preliminarily, at a specific example. Various situations of a geometrical nature have been studied successfully using operator algebra techniques in the context of Connes' noncommutative geometry~\cite{Connes}. Many such situations arise when taking quotients of topological spaces of which the quotient topology carries insufficient information or is even trivial. An example of this kind, and the one with which we shall be concerned, is that of the space of Penrose tilings of the plane (see section~\ref{sec:2}), which can be construed as a quotient of Cantor space (more precisely, of a subspace $K$ of Cantor space homeomorphic to it), giving rise to an AF C*-algebra $\pena$ whose dimension group $\mathbb Z+\tau\mathbb Z$ (where $\tau$ is the golden number $\frac{1+\sqrt 5}2$) contains information about the frequencies of appearance of finite patterns in an arbitrary tiling, and whose equivalence classes of irreducible representations can be identified with the tilings themselves --- we can say that $\pena$ \emph{classifies} the tilings in the sense that the irreducible representations of $\pena$  ``are" the tilings. 

Of course, then the quantale $\Q \pena$ also classifies the tilings, but in this paper we use the logical ideas discussed above in order to define, in section~\ref{sec:3}, a quantale $\penq$ motivated directly by the geometry of the tilings, and in sections~\ref{sec:5} and~\ref{sec:6} we show that $\penq$ also classifies them because its \emph{relational} points correspond bijectively to the tilings. An important difference between $\Q \pena$ and $\penq$ is thus that the points of the former are representations on lattices of projections of Hilbert spaces, whereas for the latter we restrict to representations on powersets. In this way $\penq$ provides a kind of ``dynamical" space that we may see metaphorically as a ``quantum space without superposition". The relation between $\penq$ and $\pena$ is currently being studied, but will not be addressed here.

This paper also provides insights into the nature of the notion of point of an involutive quantale, because on the one hand it contains results of independent interest concerning relational representations, in particular a decomposition theorem showing that any relational representation is uniquely partitioned into irreducible components (section~\ref{sec:4}), while on the other it allows us to see, in the setting of relational representations and in a very explicit example, the way in which algebraically irreducible representations differ from those representations that are only irreducible. Although it has been proposed elsewhere~\cite{PeRo97} (see also~\cite{Kr02}) that these latter should also be considered as giving rise to points of an involutive quantale, we see in this context that they fail to yield the intended classification of tilings, just as in the case of the spectrum $\Q A$ of a C*-algebra $A$ the intended points of the spectrum, namely the algebraically irreducible representations on Hilbert space, seem to require consideration of the algebraically irreducible representations of the quantale $\Q A$. 

In the present case, this classification arises very explicitly from a relational representation of $\penq$ on Cantor space (rather, on its subset $K$ mentioned above), of which the irreducible components are necessarily algebraically irreducible and coincide precisely with the equivalence classes obtained when the set of tilings is viewed as a quotient of $K$. Hence, in particular, the identification of the set of Penrose tilings with a quotient of $K$ is intrinsically contained in the representation theory of $\penq$.

We conclude this introduction by recalling some basic definitions:

By a \emph{sup-lattice} $S$ is
meant a partially ordered set, each of whose subsets $X\subseteq S$
has a join $\V X\in S$. Hence, in particular, a sup-lattice is in fact a complete lattice.

By a \emph{sup-lattice homomorphism} is meant a mapping $h:S\rightarrow S'$ between sup-lattices
that preserves all joins:
\[h(\V X)=\V\{h(x)\st x\in
X\}\]
for all $X\subseteq S$.

By a \emph{quantale} $Q$ is meant a sup-lattice,
together with an
associative product, $\then$,  satisfying
\[a \then \biggl(\V_i b_{i}\biggr) = \V_i (a \then b_{i}) \]
and
\[\biggl(\V_i a_{i}\biggr) \then b = \V_i (a_{i} \then b) \]
for all $a, b, a_{i}, b_{i} \in Q$. The bottom element
$\V\emptyset$ of $Q$ is denoted by
$0$, and the top element $\V Q$ is denoted by $1$. The quantale
$Q$ is said to be \emph{unital} provided
that there exists an element $e \in Q$ for which
\[e \then a = a = a \then e \]
for all $a \in Q$.

By an \emph{involutive quantale} is meant a quantale $Q$ together with an
involution, $^*$,
satisfying the conditions that
\begin{eqnarray*}
a^{**} &=& a\;,\\
(a \then b)^* &=& b^* \then a^*\;,\\
\biggl(\V_i a_{i}\biggr)^* &=& \V_i a_{i}^*\;
\end{eqnarray*}
for all $a, b, a_{i} \in Q$. An element $a\in Q$ for which $a^{*}=a$ is said to be 
\emph{self-adjoint}.

By a \emph{homomorphism} of quantales $h:Q\rightarrow Q'$ is meant a
sup-lattice homomorphism that also preserves multiplication,
\[h(a\then
b)=h(a)\then h(b)\]
for all $a,b\in Q$. The homomorphism is said to be \emph{strong} provided that
\[h(1_{Q})=1_{Q'},\] to be \emph{unital}
provided that \[h(e_{Q})\ge e_{Q'},\] and to be \emph{strictly unital} provided that \[h(e_{Q})= e_{Q'}.\]

A homomorphism $h:Q\rightarrow Q'$ of involutive quantales is said to be
\emph{involutive} provided that
$h(a^*)=h(a)^{*}$ for all $a\in Q$. 
In this paper all quantales and homomorphisms are assumed to be involutive
and unital.

\section{Penrose tilings}\label{sec:2}

We begin by recalling the facts that we shall need about Penrose tilings.
Generically, a tiling of the plane (which for convenience we 
identify with the complex plane) is a covering
\[\mathbb{T}=(T_i)_{i\in I}\]
of
the plane by connected closed subsets, satisfying suitable 
conditions~\cite{T&P}. In the case of Penrose tilings we shall make those conditions
explicit in a moment.
It may be remarked that we may define, on the set of all the tilings of the plane,
an action of the (additive group of) complex numbers by translation: for
each tiling $\mathbb T$ and each $z\in \mathbb{C}$ define $\mathbb
T+z$ to be the family $(T_i+z)_i$, where $T_i+z$ is given by pointwise
addition $\{w+z\in\mathbb C\st w\in T_i\}$.

The Penrose tilings that we shall be considering, which form a set that
is closed under the action by translations just described, use only two
basic shapes: namely, the two triangles depicted in Fig.~\ref{fig:basictiles}.

\setlength{\unitlength}{1mm}
\thicklines

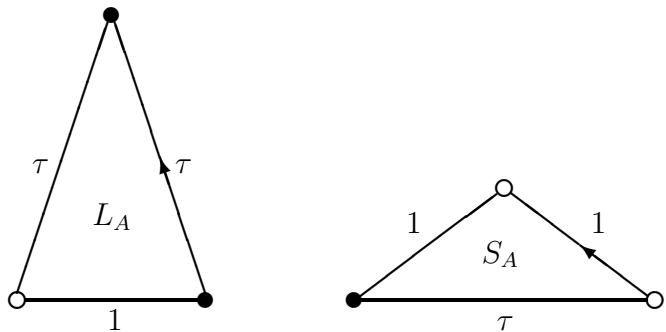
\begin{figure}
\begin{center}
\begin{picture}(27,50)
\put(0,0){\circle{2}}
\put(25,0){\circle*{2}}
\put(12.5,37.9){\circle*{2}}
\put(1,0){\line(1,0){23}}
\put(0,1){\line(1,3){12}}
\put(25,1){\vector(-1,3){6}}
\put(19,19){\line(-1,3){6}}
\put(10,10){$L_A$}
\put(12,-4){$1$}
\put(2,17){$\tau$}
\put(21,17){$\tau$}
\end{picture}
\hspace*{1.5cm}
\begin{picture}(41,17)
\put(0,0){\circle*{2}}
\put(40,0){\circle{2}}
\put(20,15){\circle{2}}
\put(1,0){\line(1,0){38}}
\put(0.9,0.6){\line(4,3){18.1}}
\put(39.1,0.6){\vector(-4,3){9}}
\put(29.9,7.47){\line(-4,3){9}}
\put(17,5){$S_A$}
\put(19,-4){$\tau$}
\put(7,9){$1$}
\put(31.5,9){$1$}
\end{picture}
\end{center}
\caption{Basic tiles of a Penrose tiling ($\tau=\frac{1+\sqrt 5}2$).}\label{fig:basictiles}
\end{figure}

The vertices of these triangles are coloured and some edges are oriented
as shown, and in a tiling both the colours of the coinciding vertices and
the orientations of the coinciding edges must agree. The triangle $L_A$
is the \emph{large tile} and $S_A$ is the \emph{small tile}, where large
and small refer to their respective areas. It is known that the whole
plane can be tiled in this way~(\cf~\cite{T&P}).

\begin{figure}
\begin{center}
\setlength{\unitlength}{1.3mm}
\begin{picture}(41,17)
\put(0,0){\circle*{1.3}}
\put(40,0){\circle{1.3}}
\put(20,15){\circle*{1.3}}
\put(0.7,0){\line(1,0){22.3}}
\put(39.3,0){\vector(-1,0){8}}
\put(31.3,0){\line(-1,0){7}}
\put(0.35,0.25){\line(4,3){9.45}}
\put(19.65,14.7){\vector(-4,-3){9.75}}
\put(39.4,0.3){\line(-4,3){19}}
\put(20,15){\line(1,-4){3.55}}
\put(23.8,0){\circle{1.3}}
\end{picture}
\hspace*{1cm}
\begin{picture}(41,17)
\put(0,0){\circle*{1.3}}
\put(40,0){\circle{1.3}}
\put(20,15){\circle*{1.3}}
\put(0.7,0){\line(1,0){38.6}}
\put(0.35,0.25){\line(4,3){9.45}}
\put(19.65,14.7){\vector(-4,-3){9.75}}
\put(39.4,0.3){\line(-4,3){19}}
\put(17,5){$L_B$}
\put(14,-4){$\tau+1=\tau^2$}
\put(7,9){$\tau$}
\put(31.5,9){$\tau$}
\end{picture}
\end{center}
\vspace*{1cm}
\begin{center}
\setlength{\unitlength}{0.8mm}
\begin{picture}(27,37)
\put(0,0){\circle{2}}
\put(25,0){\circle*{2}}
\put(12.5,37.9){\circle*{2}}
\put(1,0){\line(1,0){23}}
\put(0,1){\line(1,3){12}}
\put(25,1){\vector(-1,3){6}}
\put(19,19){\line(-1,3){6}}
\put(10,10){$S_B$}
\put(28,10){$= L_A$}
\put(12,-4){$1$}
\put(2,17){$\tau$}
\put(21,17){$\tau$}
\end{picture}
\end{center}
\caption{Constructing the two type $B$ tiles.}\label{fig:typeBtiles}
\end{figure}
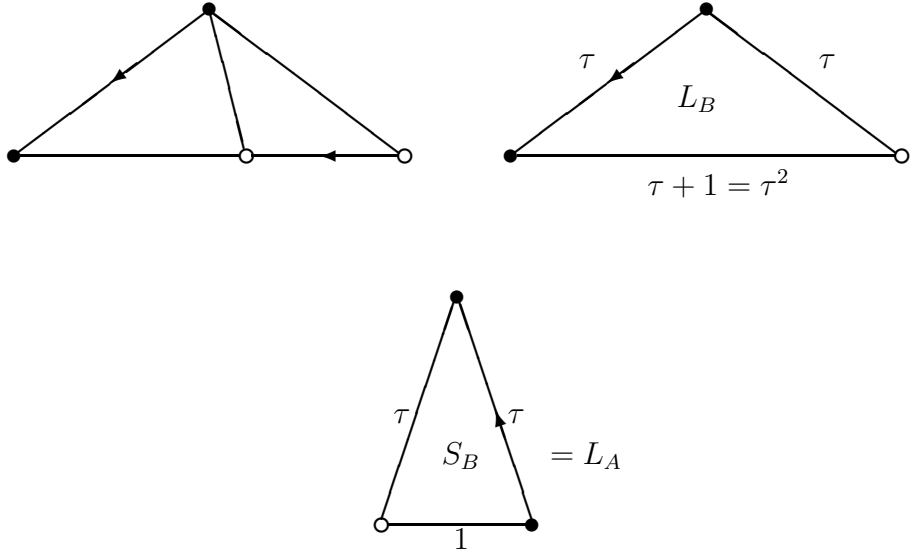

A feature of this kind of tiling is that
they possess a certain invariance of scale whereby each tiling determines
a denumerable family of other tilings, as we now explain. First, notice
that for the configuration of Fig.~\ref{fig:typeBtiles} we can, by removing the edge
between the small and the large tile, produce a new tile congruent to
$S_A$, with the colours of the vertices and the orientations shown.

\begin{figure}
\begin{center}
\includegraphics[width=0.45\textwidth]{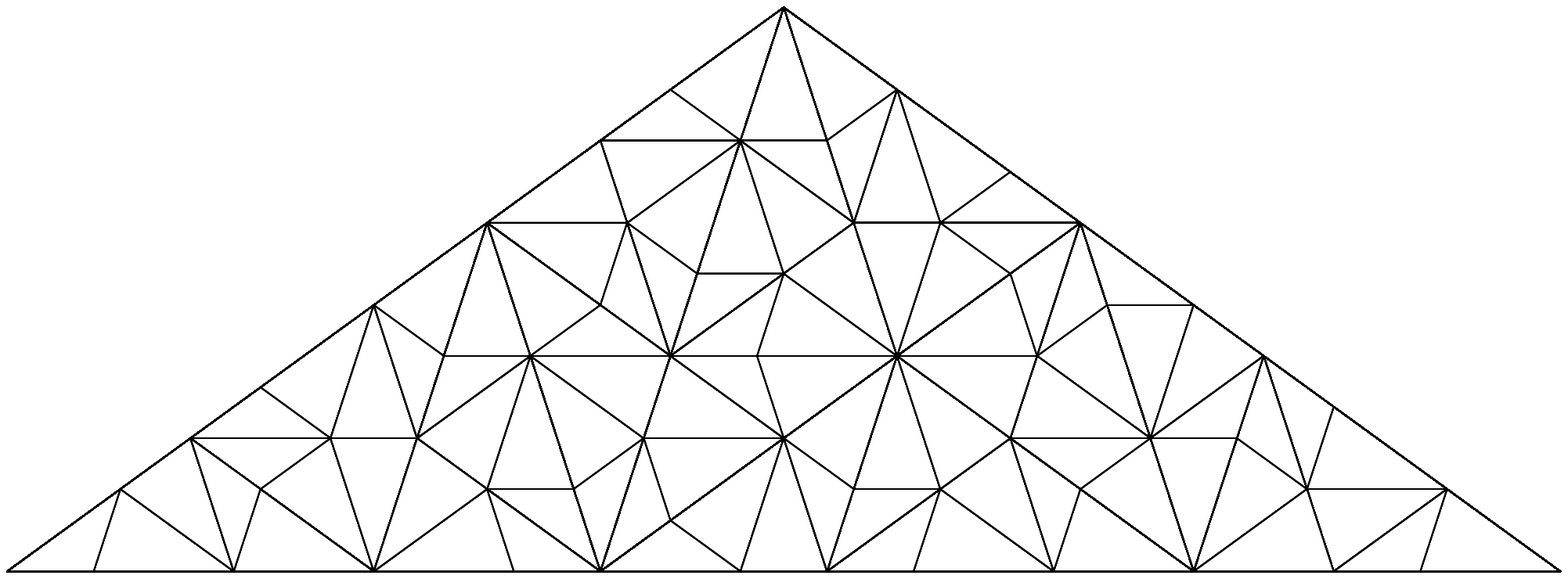}
\includegraphics[width=0.45\textwidth]{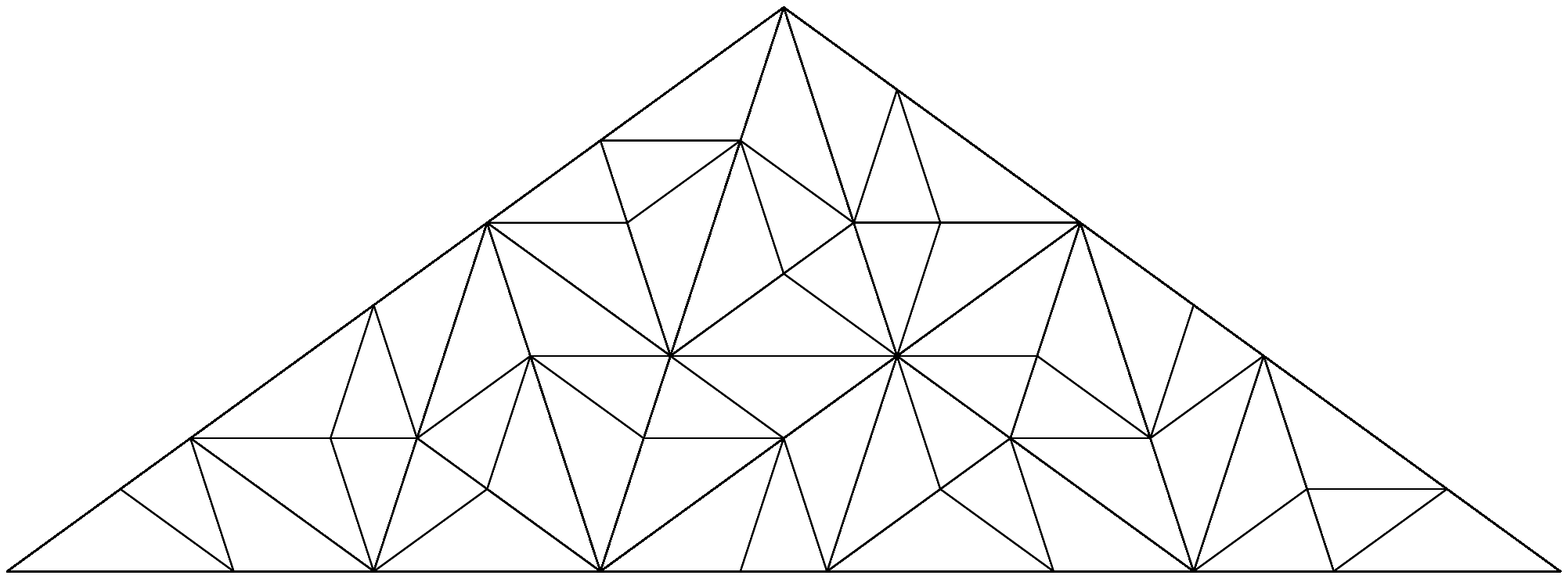}
\end{center}
\vspace*{-2cm}
\caption{Left: Finite fragment of a Penrose tiling. Right: Same fragment after edges were deleted so as to obtain a tiling by type $B$ tiles.}\label{fig:tilings}
\end{figure}

It is now possible, following the same rules as before, to tile the whole
plane using this tile and $L_A$, which are now named $L_B$ and $S_B$,
respectively, because $L_A$ has become the smaller tile. The subscript
$B$ refers to the new type of tiling, which differs from the previous one
(type $A$) both due to the relative sizes of the two triangles, and 
to the different arrangements in the colours of their vertices and the
orientations of their edges. Fig.~\ref{fig:tilings} shows both a fragment of a tiling of type $A$ and the same fragment after edges were removed so as to obtain a tiling of type $B$.

There is another important aspect to bear in
mind, namely that each tile of the original tiling is contained in a
unique tile of the new tiling, which we refer to as its \emph{successor}.
More than that, the successor of a large tile can be large or small, but
the successor of a small tile must necessarily be large.

In the same way that we delete certain edges between tiles of type $A$ in order to obtain those of type $B$, we may delete certain edges between tiles of type $B$ in order to obtain larger tiles and yet another tiling of the plane. Furthermore, this procedure can be
continued indefinitely, yielding a sequence of tilings
$(\mathbb{T}_n)_{n\in\N}$, where each $\mathbb{T}_n$ is said to be a tiling of
\emph{level} $n$ (levels 0 and 1 correspond to the tilings of type $A$
and $B$, respectively). Similarly to the passage from type $A$ to type
$B$, each tile of level $n$ must be contained in a unique tile of level
$n+1$, again called its \emph{successor}, and the successor of a small
tile must necessarily be a large tile, whereas the successor of a large
tile can be either small or large; it is small precisely if it coincides
with the large tile of level $n$. In order to make our notation more
concise we shall denote the successor of a tile $T$ by $s(T)$. Successors commute with translations,
\ie, $s(T+z)=s(T)+z$ for any
tile $T$ and any $z\in\mathbb C$.

Given such a tiling and a point $z$ on the complex plane located in 
the interior of a tile $T$ (at level zero), we define a 
sequence $\tilde z$ of zeros and ones as follows:
\[\tilde z_{n}=\left\{\begin{array}{ll}
0 & \textrm{if }s^{n}(T)\textrm{ is large at level }n\\
1 & \textrm{if }s^{n}(T)\textrm{ is small at level }n\;.
\end{array}\right.
\]
Since every small tile at a level $n$ becomes part of a large tile at 
level $n+1$ we conclude that the sequence $\tilde z$ has the property
that every $1$ must be followed by a 
$0$.

\begin{definition}\label{def:Penrose.sequence}
    By a \emph{Penrose sequence} is meant a sequence $s\in\{0,1\}^\N$ satisfying the condition that
    \[s_n=1\textrm{ implies that }s_{n+1}=0\]
    for all $n\in\N$. The set of all Penrose sequences is denoted by $K$.
\end{definition}

It may be shown that any Penrose sequence can be obtained as described 
above from a Penrose tiling, and that two sequences $s$ and $t$ are
obtained from two points 
on the same tiling if, and only, if they are equivalent in the 
following sense~\cite{T&P}:

\begin{definition}\label{def:equivalence.sequences}
    The Penrose sequences $s$ and $t$ are said to be \emph{equivalent}, 
    written $s\sim t$, provided that
    for some $m\in\mathbb N$ one has that $s_{n}=t_{n}$ for all $n\ge m$. 
    (Equivalently, $s$ and $t$ differ only in finitely many places.)
\end{definition}

These facts lead to the identification of the set of Penrose tilings 
with the quotient set $K/{\sim}$ of the set $K$ determined by this equivalence relation. In fact, $K$ is 
naturally equipped with a topology that is generated by subbasic open sets of 
the form
    \[U(n,b)=\{s\in K\st s_{n}=b\}\]
for each $n\in\mathbb 
N$ and $b\in\{0,1\}$. This is the subspace topology obtained from the 
Cantor topology by identifying $\{0,1\}^{\mathbb 
N}$ with Cantor space, and indeed $K$ is homeomorphic to Cantor space.
However, considering
$K$ as a topological space rather than just a set is of no use from 
the point of view of 
Penrose tilings, because all the equivalence classes are dense and 
hence the quotient topology is trivial. This provides one of the motivations 
for considering the set of Penrose tilings instead as a generalised kind of
topological space, for instance by identifying it with a
C*-algebra that plays the r\^{o}le of a ``noncommutative space", as in~\cite{Connes}.

\section{A theory of Penrose tilings}\label{sec:3}

In order to treat the set of Penrose tilings as a quantale, let us
present one by means of propositions and axioms of a logical
theory which is noncommutative in the sense that conjunction is a
noncommutative connective. To a great extent this section will be devoted to
providing some level of intuitive motivation for the axioms we choose, and as such it is not entirely mathematical. The reader who does not wish to be distracted by such considerations may read the definition of the theory below, and move directly to section~\ref{sec:4}. The results in sections~\ref{sec:5} and~\ref{sec:6} can also be regarded as an \textit{a posteriori} justification for the axioms, at which point the motivational reasoning takes on a more mathematical significance.

The propositions of the theory with which we are concerned should be
viewed as representing ``experiments", or ``measurements", by means
of which we as observers learn about the structure of a particular, but arbitrary,
Penrose tiling. Furthermore, the tiling being studied remains fixed, up to
translations, during the whole series of experiments.
Each proposition is built by applying appropriate connectives to primitive propositions. We list the primitive propositions below, followed in each case by the description of the experiment that it represents. In each
case we denote by ${}^n{T}$ the tile of level $n$ that contains the origin of
the plane. Throughout we restrict to tilings for which $0$ does not lie
on an edge or a vertex of any tile of level zero (hence, of any tile of any level):

\begin{itemize}
\item $\wprop n L$ --- perform a random translation $z\in\mathbb C$
 (or sequence of 
translations whose sum is $z\in\mathbb C$) such
that in the end the origin of the plane remains in the successor of ${}^n{T}$, \ie, $0\in
s({}^n{T}+z)$,
without ending, as stressed above, on an edge or a vertex of any tile; the experiment succeeds if the
origin is found inside a large tile of level $n$ after the translation
takes place.
\item $\wprop n S$ --- analogously, but with the origin being found
inside a small tile after the translation; notice that this experiment
never succeeds in those cases where ${}^n{T}$ is large and $s({}^n{T})={}^n{T}$.
\item $\rprop n L$ --- this experiment is the dual (in the sense of temporal reversal)
of $\wprop n L$; it succeeds if ${}^n{T}$ is a large tile at level $n$, in
which case the tiling is translated randomly, with the only restriction
being that the origin must remain inside $s({}^n{T})$ (and, as before, not on an edge or vertex).
\item $\rprop n S$ --- similar to the previous case, but with ${}^n{T}$ being small.
\end{itemize}
It should be noted that we have here, in order to make their interpretation explicitly available, considered both $\wprop n L$
and $\rprop n L$ as primitive propositions, and similarly for $\wprop n S$ and $\rprop n S$, although the axioms that we shall now introduce identify these as being obtained each from the other by applying the logical operator of temporal dual introduced below. Of course, the theory may equivalently be presented just in terms of the primitive propositions $\wprop n L$ and $\wprop n S$.

Based on this geometric interpretation of the primitive propositions
we may now elicit
some axioms relating them. As stated above we shall use a conjunctive 
connective, $\then$, meaning
``and then",
where \[a\then b\] is the experiment that
consists of performing first $a$ and then $b$, together with a disjunctive connective, $\V$, which can be infinitary (or, more precisely, a family of disjunctive connectives, one for each arity), where
\[\V_{i\in I} a_i\] is the experiment that consists of performing $a_i$ for some $i\in I$.
Finally, we introduce a ``temporal 
dual'' \[a^{*}\] for each proposition $a$, with the dual of $\wprop n L$ being $\rprop n L$, with
$\wprop n S$ similarly having dual $\rprop n S$. For an introduction to the logical background for this noncommutative propositional logic the reader is referred to~\cite{Mu86}, and for its application in another context to~\cite{MuPe92,Re02:jpaa}. A similar context, but without the temporal dual operator, may be found in~\cite{AbVi93,Re99,Re00,Re01,ReVi}.

The logic also introduces the truth value $\false$, here representing the impossible experiment that always fails, together with the truth value $\true$, corresponding here to the experiment that always succeeds; more precisely, the experiment that either produces no translation at all or which, having produced a translation, leaves the tiling back in its original position and provides no information whatsoever about the type of tiles in which the origin is located.

The axioms of the theory will be written in the form
\[a\entails b\;,\]
meaning that ``$a$ implies $b$": in this context, that the experiment
$a$ is a particular way of performing the (more general) experiment $b$.
More precisely, this means that if $a$ succeeds then $b$ succeeds also,
and whatever information is extracted from $a$ could also have been obtained by
performing $b$.

We now provide a list of axioms that gives an intuitive description of the way in which
measurements represented by the primitive propositions relate to each other. In section~\ref{sec:5} we shall see that,
in an appropriate sense, these axioms provide a complete description of the geometric situation. Throughout, as abbreviations we
write:
\begin{itemize}
\item $X$ or $Y$ for either $S$ or $L$;
\item $\langle \bprop n X\vert\bprop m Y\rangle$ instead of $\wprop n X\then\rprop m
Y$;
\item $\vert \bprop n X\rangle\langle \bprop m Y\vert$ instead of
$\rprop n X\then\wprop m Y$;
\item $\langle \bprop{n_{1}}{X_{1}}; \ldots; \bprop{n_m}{X_{m}}\vert$ instead of
$\wprop{n_1}{X_{1}}\then\ldots\then\wprop{n_m}{X_{m}}$;
\item $\vert \bprop{n_{1}}{X_{1}}; \ldots; \bprop{n_{m}}{X_{m}}\rangle$ instead of
$\rprop{n_1}{X_{1}}\then\ldots\then\rprop{n_m}{X_{m}}$.
\end{itemize}
We shall also say that a string $\bprop 0{X_{0}};\ldots; \bprop 
n{X_{n}}$ is
\emph{admissible} if $X_{i}=S$ entails 
$X_{i+1}=L$, for all $i\in\{0,\ldots,n-1\}$.

\begin{definition}\label{def:pent}
The \emph{theory $\pent$ of Penrose tilings} is given by taking the primitive propositions
\[\wprop n L\textrm{ and }\wprop n S\textrm{ (together with their duals }\rprop n L\textrm{ and }\rprop n S\textrm{)}\]
for each $n\in\N$, together with the following axioms:
\[
\begin{array}{lll}
(C1_{n}) &\langle \bprop n L\vert \bprop n S\rangle\entails
\false&\textrm{(Consistency 1)}\\
(C2_{n}) &\wprop n S\entails\wprop {n+1}L&
\textrm{(Consistency 2)}\\
(D1_{n}) &\true\entails\wprop n L\vee\wprop n S&
\textrm{(Decidability 1)}\\
(D2_{n}) &\vert \bprop{n+1}X\rangle\langle \bprop{n+1}X\vert\entails
\wprop n S\vee\wprop n L&
\textrm{(Decidability 2)}\\
(E1_{n}) &\langle \bprop n Y; \bprop{n+1}X\vert\entails \wprop {n+1}X&
\textrm{(Expansion 1)}\\
(E2_{n}) &\langle\bprop{n+1}X;\bprop n Y\vert\entails \wprop 
{n+1}X&\textrm{(Expansion 2)}\\
(E3_{n}) &\rprop n Y\wprop{n+1}X\entails \wprop {n+1}X&
\textrm{(Expansion 3)}\\
(E4_{n}) &\langle\bprop{n+1}X\vert\bprop n Y\rangle\entails \wprop 
{n+1}X&\textrm{(Expansion 4)}\\
(I1_{n}) &\rprop n X^*\entails\wprop n X&
\textrm{(Involution 1)}\\
(I2_{n}) &\wprop n X^*\entails\rprop n X&
\textrm{(Involution 2)}\\
(C'_{t}) &\true\entails \langle \bprop n{{X}_n};\ldots; \bprop 0{{X}_0}\vert
\bprop 0{{X}_0};\ldots;\bprop n{{X}_n}\rangle&
\textrm{(Completeness)}
\end{array}
\]
for all 
$n\in\N$, all
$X,Y\in\{S,L\}$, and all admissible strings 
$t=\bprop 0{X_{0}};\ldots; \bprop 
n{X_{n}}$ with $X_{n}=L$. Again, it should be noted that the axioms $(I1_n)$ and $(I2_n)$ may be omitted if the primitive propositions $\rprop n L$ and $\rprop n S$ are considered to have been added formally by the logical context.
\end{definition}

We shall now justify each of these axioms in terms of the geometric 
interpretation of the primitive propositions:

\begin{description}
    \item[\textnormal{($C1_{n}$):}] The proposition $\langle \bprop n L\vert 
    \bprop n S\rangle$ consists of performing $\wprop n L$ and then 
    $\rprop n S$. But $\wprop n L$ fails if the tile of level 
    $n$ where the 
    origin is located when the experiment ends is small, whereas 
    $\rprop n S$ fails unless it starts with the origin located in a 
    small tile at level $n$. Hence, the sequence of the two experiments always fails, 
    \ie, it entails the proposition $\false$.

    \item[\textnormal{($C2_{n}$):}] The proposition $\wprop n S$ succeeds if
    a translation of the tiling that leaves the 
    origin inside its original tile of level $n+1$ is performed, such that in 
    the end the origin is found in a small tile at level $n$. This implies 
    that in the end it is also located in a large tile of level $n+1$, and thus 
    $\wprop n S$ is a particular way of doing the experiment 
    $\wprop{n+1}L$ because on one hand the translation it involves 
    certainly leaves the origin inside its tile of level $n+2$, and 
    when the experiment finishes the origin is found in a large tile 
    at level $n+1$.

    \item[\textnormal{($D1_{n}$):}] The experiment $\true$ always succeeds, and for 
    an arbitrary $n$ it is either a particular way of doing $\wprop n S$ or of 
    doing $\wprop n L$, because the origin must always be either in a small 
    tile or in a large tile at level $n$.

    \item[\textnormal{($D2_{n}$):}] The experiment
    $\vert \bprop{n+1}X\rangle\langle \bprop {n+1}X\vert$ consists of doing
    $\rprop {n+1}X$ and then $\wprop{n+1}X$. In either of these two 
    experiments
    the translations involved leave the origin inside its triangle at 
    level $n+2$, and the whole experiment succeeds if, and 
    only if, the tile at level $n+1$ in which the origin initially 
    lies is of the same type (large or small) as that where it lies 
    when the experiment ends. Since inside a tile of level $n+2$ 
    there is at most one tile of each type at level $n+1$, we 
    conclude that the translation involved in performing the whole 
    experiment leaves the origin inside its tile of level $n+1$. 
    Hence, if the origin ends up located in a large tile at level $n$ 
    we will have performed the experiment $\wprop n L$, otherwise we 
    will have performed $\wprop n S$, which means that
    $\vert \bprop {n+1}X\rangle\langle \bprop {n+1}X\vert$ is a way of performing 
    the disjunction of the two.

    \item[\textnormal{($E1_{n}$):}] The experiment $\langle \bprop n Y; \bprop {n+1}X\vert$
    consists of performing $\wprop n Y$ and then $\wprop{n+1}X$, which 
    means first performing a translation that leaves the origin in its 
    triangle of level $n+1$, inside a triangle of type $Y$ at level 
    $n$, and then performing a translation that leaves the origin in 
    its triangle of level $n+2$ and in a triangle of type $X$ at 
    level $n+1$. This sequence, of course, provides a particular way
    of performing $\wprop{n+1}X$.

    \item[\textnormal{($E2_{n}$):}] Similarly, but $\wprop {n+1}X$ is performed first.

    \item[\textnormal{($E3_{n}$):}] Similar to $E1_{n}$, except that the first 
    translation is performed only if the origin of the plane is 
    initially inside a triangle of type $Y$ at level $n$, and after 
    the first translation it can be in either of the two types at 
    level $n$.

    \item[\textnormal{($E4_{n}$):}] Similar, but $\wprop {n+1}X$ is performed first.

    \item[\textnormal{($I1_{n}$) and ($I2_{n}$):}] These axioms express that $\wprop n 
    X$ is the time reversal of $\rprop n X$.

    \item[\textnormal{($C'_{t}$):}] In order to understand this axiom, let us first 
    see what the effect of the experiment $\langle\bprop n{{X}_n};\ldots;
    \bprop 0{{X}_0}\vert$ is. Recall that we are assuming that $X_{n}=L$.
    First, $\wprop n {X_{n}}$, \ie, $\wprop n {L}$, will translate the tiling and 
    leave the origin inside its original tile of level $n+1$, and in a 
    large tile at level $n$ (notice that this experiment always 
    succeeds --- it might fail if $X_{n}$ were $S$ instead, which is 
    why
    we have ruled out this possibility in the definition of the axiom);
    then, $\wprop {n-1} {X_{n-1}}$ translates again 
    the tiling, without taking the origin away from the tile of 
    level $n$ that resulted from the previous experiment, and further 
    leaving 
    it in a tile of type $X_{n-1}$ at level $n-1$; this experiment 
    succeeds for any value of $X_{n-1}$ because $X_{n}=L$. Proceeding in 
    this way it is clear
    that each of the steps in the sequence succeeds (because the 
    sequence is admissible) and that 
    the
    whole sequence succeeds if, and only if, it is possible to produce a 
    translation that does not place the origin outside its tile of 
    level $n+1$, in the end placing it in a tile of 
    type $X_{0}$ at level $0$, which in turn is inside a tile of type 
    $X_{1}$ 
    at level $1$, etc. The sequence
    $\vert {X}_0;\ldots; 
    {X}_n\rangle$ is then the time reversal of this; it succeeds if 
    and only if the origin is initially in a ``tower'' of tiles precisely equal to 
    the one just described. Hence, after $\langle\bprop n{{X}_n};\ldots;
    \bprop 0{{X}_0}\vert$ it necessarily succeeds,
    in the end producing a translation that will 
    keep the origin in the original tile of level $n+1$, and allowing us to conclude two things: first, this experiment gives us no
information whatsoever about where we are in a tiling (\ie, in which kind
of triangle at each level), either before or after the experiment takes
place; second, a particular way of performing it consists of doing a
translation that places the tiling back in its original position. But
the experiment that does such a translation and provides no information about the kind of tiles
where the origin is located
is what we defined to be the primitive proposition $\true$, which is thus
seen to be a particular way of performing $\langle \bprop n{{X}_n};\ldots; \bprop 0{{X}_0}\vert
\bprop 0{{X}_0};\ldots;\bprop n{{X}_n}\rangle$.
\end{description}

By a \emph{model} of the theory $\pent$ should now be understood, following~\cite{MuPe01}, a quantale $Q$ together with an \emph{interpretation of $\pent$ in $Q$}, by which we mean an assignment to each primitive proposition $\pi$ of $\pent$ of an element $\pi_Q\in Q$, validating the axioms of the theory, where the conjunctive connective is interpreted as multiplication, disjunctions are interpreted as joins, the temporal dual operator is interpreted as involution, $\true$ is interpreted as $e$, and $\false$ as $0$. For instance, the axiom ($C1_n$) is validated by such an interpretation if, and only if, \[\wprop n L_Q\then\rprop n S_Q\le 0\;,\]and the axiom ($D1_n$) is validated if, and only if,
\[e\le\wprop n L_Q\vee\wprop n S_Q\;.\]

As usual in logic, an interpretation of the theory $\pent$ in $Q$ can be defined equivalently to be a homomorphism to $Q$ from the Lindenbaum algebra of the theory. We can think of the latter as being the quantale obtained by taking the set of propositions in the theory modulo provable equivalence in the theory, partially ordered by provable entailment in the theory, albeit taking into account that the ``set" of propositions in the theory $\pent$ is actually not a set but rather a proper class, due to the unbounded arity of disjunction, and that provable entailment must be defined in terms of possibly infinite derivations, for the same reason (see~\cite{Re99} or~\cite{Re00}). Concretely, the ``Lindenbaum quantale" of the theory $\pent$ can be constructed in the same way as any quantale which is presented by generators and relations (see, \eg,
\cite{AbVi93,Re02:jpaa}). This gives us both a quantale and an interpretation of $\pent$ in it with the appropriate universal property, namely, that any other interpretation factors through the universal one and a unique quantale homomorphism. (Indeed, we may take this universal property as the definition of the Lindenbaum quantale.) With these provisos in mind we define:

\begin{definition}\label{def:penq}
By the \emph{quantale of Penrose tilings} will be meant the involutive unital quantale \[\penq\] obtained by taking the Lindenbaum quantale of the theory $\pent$, given explicitly by taking the set of propositions in the theory modulo provable equivalence in the theory, partially ordered by provable entailment in the theory.

Equivalently, the quantale $\penq$ may be considered to be that whose presentation by generators and relations is obtained from the theory $\pent$. Explicitly, the primitive propositions
\[\wprop n X\textrm{ (and } \rprop n X\textrm{)}\]
are considered to be the generators of $\penq$, and each axiom
\[a\entails b\]
is considered to be a defining relation $a\le b$, with the logical connectives interpreted as quantale operations in the obvious way, with ``and then" as the multiplication, with ``disjunctions" as joins, and with ``time reversal" as the involution. The logical constants $\true$ and $\false$ are interpreted respectively by the multiplicative unit $e$ and the bottom element $0$ of the quantale $\penq$.
\end{definition}

Henceforth, we shall identify each primitive proposition with the element of the quantale $\penq$ which it determines, using the same symbols and conventions as in the definition of the theory given above to denote the quantale operations.

By straightforward extensions of the observations made above in motivating the theory, one has immediately the following:

\begin{proposition}\label{prop:penq}
    The following conditions are validated in the quantale $\penq$ of Penrose tilings:
    \begin{enumerate}
	\item Each of the axioms of the theory $\pent$ of Penrose tilings;
	\item $\langle \bprop n S\vert \bprop n L\rangle\entails\false$;
	\item $\true\entails\rprop n L\vee\rprop n S$;
	\item $\langle\bprop m Y;\bprop{n}X;\bprop k Z\vert\entails \wprop 
	{n}X$  whenever $m,k < n$;
	\item $\true\entails \langle \bprop n{{X}_n};\ldots; \bprop 0{{X}_0}\vert
 	\bprop n{{X}_n}\rangle$.
    \end{enumerate}
\end{proposition}

We conclude this section by discussing 
what would have happened if we had followed the same line of thinking as above while 
changing the meaning of the primitive propositions by requiring them not 
to produce any translations at all, in other words in which the experiments have only very limited interaction with the geometric situation. The aim is to show that this 
change leads naturally 
to a topological space rather than a quantale.

More precisely, we shall think of the primitive propositions exactly in the 
same way as before, but restricting the set of available translations 
to just $0$. Hence, $\rprop n L$ and $\wprop n L$ will have 
the same meaning, as measuring the size of the tile after or before 
the $0$ translation is the same. Each of these primitive 
propositions just asserts that
the origin of the plane is inside a large tile of level $n$, and 
produces no change at all. Hence, 
there is no longer any point in distinguishing between them, and we will 
just write $\prop n L$ for both. Similarly, we will write
$\prop n S$ for $\wprop n S$, whose meaning is now the same as that of 
$\rprop n S$. In general, it is clear that we may have as axioms 
the following, for any propositions $a$ and $b$ derived from the primitive 
ones by means of multiplication, involution and joins:

\begin{itemize}
    \item triviality of time reversal: $a\entails a^{*}$;
    \item $true$ is the top: $a\entails \true$ (since any experiment is now a 
    particular instance of $\true$, that is, of doing nothing to a 
    tiling, because it produces no translations);
    \item idempotence of ``and then'': $a\entails a\then a$ (if we 
    perform an experiment $a$ once then we can repeat it as many 
    times as we like).
\end{itemize}

These axioms alone tell us that the quantale presented by the 
generators $\prop n L$ and $\prop n S$, with the axioms as defining 
relations, is in fact a locale, since any idempotent quantale whose top is 
the multiplicative unit is a locale~\cite{JoTi84}. Let us now add some 
axioms pertaining to the specific meaning of the generators. The ones 
below correspond to the first three of Definition~\ref{def:pent}, and their 
justification in the present setting is obvious. Now we write 
$a\wedge b$ instead of $a\then b$ to emphasize that $\then$ is no 
longer ``and then'' but instead just ``and'':

\begin{itemize}
    \item $\prop n L\wedge \prop n S\entails\false$.
    \item $\prop n S\entails\prop {n+1}L$.
    \item $\true\entails\prop n L\vee\prop n S$.
\end{itemize}

The locale $\cantor$ presented in this way is spatial (because it has a 
presentation without infinite joins in the relations and is thus 
coherent~\cite{Johnstone}), and its points can be identified in an obvious way
with Penrose sequences: a point $p:\cantor\rightarrow 2$
corresponds to the sequence $\hat p$ defined by
\[
\hat p_{n}=\left\{\begin{array}{ll}
0 & \textrm{ if }p(\prop n L)=1,\\
1 & \textrm{ otherwise.}
\end{array}\right.
\]
The axioms ensure that these sequences are well defined.
The locale itself is isomorphic to
the Cantor topology of the set $K$ of Penrose sequences
discussed earlier in this section (the subbasic open sets $U(n,0)$ 
and $U(n,1)$
are the extensions of the propositions $\prop n L$ and $\prop n S$, 
respectively).

From the logical point of view
that we have been trying to put forward, the space of 
Penrose sequences thus corresponds to the theory of how to observe an 
arbitrary Penrose 
tiling when the origin is located inside a particular tile that remains 
fixed. In view of this, it is only natural that each tile gives rise to a 
different model of the theory, \ie, a different point of the space. With hindsight, then, we may
regard the presentation of the locale $\cantor$ by generators and relations as being a theory of Penrose sequences themselves, which is made obvious if we replace the symbols $\bprop n L$ and
$\bprop n S$ by the assertions \[(s_n=0)\textrm{ and }(s_n=1)\;,\] respectively, concerning a generic Penrose sequence $s$:

\begin{itemize}
    \item $(s_n=0)\wedge (s_n=1)\entails\false$.
    \item $(s_n=1)\entails (s_{n+1}=0)$.
    \item $\true\entails (s_n=0)\vee (s_n=1)$.
\end{itemize}

To conclude, we look again at the axioms of Definition~\ref{def:pent}. 
Besides the first three, which made their way into the presentation of 
the locale $\cantor$, all the others except the completeness axiom become 
trivial. As examples we point out the following two:
\begin{itemize}
\item the first expansion axiom just states that
$\prop {n}X\wedge\prop {n+1} Y\le\prop{n+1} Y$, which of course is true 
in any locale;
\item the second decidability axiom is equivalent to the condition
\[\prop{n+1}X\le\prop n L\vee\prop n S\;,\] which now follows from the 
first decidability axiom and the fact that $\true$ is the top, for
$\prop{n+1}X\le\true\le\prop n L\vee\prop n S$.
\end{itemize}

Hence, we see that most of the axioms are essentially ``local" in 
character.
The only exception is the completeness 
axiom, which in a locale would entail that for all primitive 
propositions $\prop n X$ we must have $\true\le\prop n X$. The locale 
presented in this way has only two elements, \ie, it is the topology 
of the space with only one point, and thus we see that the completeness
axiom is the only one that makes a truly 
non-local assertion about Penrose tilings, leading us to the 
same trivial topology that one obtains when trying to describe 
the space of Penrose tilings as a quotient of Cantor space.

Another way of expressing this is to observe that consideration of the natural equivalence that arises between Penrose tilings under symmetries of the plane is intrinsically incompatible with the concept of topological space in its usual commutative sense, or indeed the concept of locale which is its constructive counterpart. To allow its expression, one needs to work with noncommutative spaces, in this context by bringing in the concept of quantale, hence equivalently by working within the context of noncommutative logic.

\section{Representation theory of quantales}\label{sec:4}

In this section we first recall some basic facts and definitions about representations of quantales, and then study relational representations, which are the ones that we shall need in this paper. We remark that relational representations are well behaved, in the sense that they can always be decomposed into irreducible components (Theorem~\ref{thm:irredreps}), in such a way that in order to describe completely the relational representations of any quantale it suffices to know those that are irreducible.
This situation does not hold in general for other kinds of representations of quantales.

Let $S$ be an orthocomplemented sup-lattice, by which we mean a sup-lattice 
equipped with an antitone automorphism $x\mapsto x^\perp$
satisfying the following two conditions:
\begin{eqnarray*}
    x &=& x^{\perp\perp}\;,\\
    0 &=& x\wedge x^\perp\;.
\end{eqnarray*}
We write $x\perp y$, and say that $x$ and $y$ are 
\emph{orthogonal}, if $x\le y^\perp$ (equivalently, $y\le x^\perp$). 
As examples of this we note in particular the sup-lattice $\pwset H$ of projections of a 
Hilbert space $H$, with the usual orthogonal complement, and the 
powerset $\pwset X$
of any set $X$, with the set-theoretic complement $Y\mapsto X\setminus Y$.

The set $\End(S)$ of all the sup-lattice 
endomorphisms of $S$ is a quantale under
the pointwise ordering, with
multiplication given by composition, and the involution 
defined adjointly by
\[f^{*}(y)^\perp=\V\{x\in S\st f(x)\le y^\perp\}\;.\]
The reader is referred to~\cite{MuPe92,MuPe01,MuPe02,PeRo97,Re:2-forms} for background discussion of quantales of this kind.

By a \emph{representation} 
of a quantale $Q$ on the orthocomplemented 
sup-lattice $S$ is meant a homomorphism
\[r:Q\rightarrow\End(S)\;.\]
Choosing the multiplication in $\End(S)$ to be
$f\then g=g\circ f$, from the representation $r$ we obtain,
for each $x\in S$ and $a\in Q$, an action of $a$ on $x$
given by $xa=r(a)(x)$.

In this way, we obtain from any representation of a quantale $Q$ on an orthocomplemented sup-lattice $S$, a
 \emph{right action} of $Q$ on $S$ in the usual sense of a semigroup acting 
on a set, which moreover preserves both the joins of $Q$ and those 
of $S$:
\begin{eqnarray*}
    (\V_{i}x_{i})a&=&\V_{i}(x_{i}a)\;,\\
    x(\V_{i}a_{i})&=&\V_{i}(xa_{i})\;,
\end{eqnarray*}
thereby making $S$ a \emph{right $Q$-module}. This module is 
\emph{unital}, meaning that $xe\ge x$ 
for all $x\in S$, since the homomorphism defining the representation is assumed unital. It is evidently \emph{strictly unital}, in the sense that $xe=x$ for all $x\in S$, precisely if the representation is 
strictly unital.

The fact that $r$ is an involutive homomorphism can be restated 
equivalently by saying that $S$ is an \emph{involutive right $Q$-module}~\cite{Re:2-forms}, 
\ie, a right $Q$-module satisfying the 
condition
\[xa\perp y\iff x\perp ya^{*}\]
 for all $x,y\in S$ and all $a\in Q$.

By a \emph{homomorphism} of right $Q$-modules $f:S\rightarrow S'$
is meant a sup-lattice homomorphism 
that commutes with the action:
 \[f(xa)=f(x)a\]
 for each $x\in S$ and $a\in Q$.
For further facts concerning involutive quantale modules and their 
homomorphisms see~\cite{Re:2-forms}.

Given any representation $r:Q\rightarrow\End(S)$, by an \emph{invariant element} of the representation is meant an element 
$x\in S$ with the property that \[x 1\le x.\]
Evidently this is equivalent to requiring that $xa\le x$ for each $a\in Q$. It is also equivalent in the present
 context, since the representation is a unital homomorphism, to requiring simply that
 \[x1=x\;.\]
 The representation is said to be \emph{irreducible} if $S$ is non-zero and there are no invariant elements besides $0$ and $1$.
 The latter condition is equivalent by a straightforward argument to the representation being \emph{strong}, that is, that
$r(1_{Q})=1_{\End(S)}$. The invariant elements of $Q$, when 
$Q$ is viewed as a right module over itself, are evidently the 
\emph{right-sided} elements of $Q$.

By a \emph{cyclic generator} of a 
representation $r:Q\rightarrow\End(S)$ is meant an element $y\in S$ for which one has that
 \[\{ya\in S\st a\in Q\}=S\;.\]
In the case when $S$ is non-zero and atomic as a sup-lattice, the representation 
is said to be \emph{algebraically irreducible}
provided that each of its atoms is a cyclic generator~\cite{MuPe01}. Any algebraically 
irreducible representation of a quantale is necessarily irreducible.

By an \emph{equivalence} of representations from a representation $r:Q\rightarrow\End(S)$
to $r':Q\rightarrow \End(S')$ is meant an isomorphism of right $Q$-modules
$f:S\rightarrow S'$ which is an ``isometry'' in the sense that it preserves and 
reflects orthogonality: \[x\perp y\iff f(x)\perp f(y)\]for all $x,y\in S$.
The representations are said to be \emph{equivalent} if there is an equivalence 
between them.

The importance of the notion of algebraic irreducibility lies in the 
observation~\cite{MuPe01} that for the quantale $\Q A$ associated to a unital 
C*-algebra $A$ there is a bijective correspondence, up to 
equivalence of representations, between the 
irreducible representations of $A$ and the algebraically irreducible 
representations of $\Q A$, which are thus regarded as 
the ``points'' of $\Q A$. It may be noted that at present this correspondence holds only up to a 
conjecture, equivalent to asking~\cite{Re:2-forms} that for each algebraically irreducible representation
of $\Q A$ the annihilator
 \[\ann(x)=\V\{a\in\Q A\st xa=0\}\]
  of at 
least one cyclic generator $x\in S$ is a maximal right-sided element of 
$\Q A$. The precise form of the conjecture in~\cite{MuPe01} is in terms of the non-triviality of the representation
  with respect to the mapping of pure states.

 In the present context, we shall be interested in the particular case of representations on 
orthocomplemented sup-lattices 
of the form $\pwset X$ for a set $X$. As we shall see, these 
representations are well behaved in the sense that
they can always be decomposed in a unique way into irreducible components.
We begin by remarking that each sup-lattice 
endomorphism \[f:\pwset X\rightarrow \pwset X\] determines a
binary relation $R_{f}$ on $X$ defined by
\[x R_{f} y\iff y\in f(\{x\})\] for each $x,y\in X$.
Conversely, each binary relation $R$ on $X$ determines a sup-lattice 
endomorphism $f_{R}:\pwset X\rightarrow\pwset X$ by defining
\[f_{R}(Z)=\bigcup_{z\in Z}\{x\in X\st z R x\}\] for each element $Z\in\pwset X$, giving
a bijective mapping from the quantale $\End(\pwset X)$ to the sup-lattice $\pwset{X\times X}$.
It is verified straightforwardly that with respect to the natural structure of quantale on the sup-lattice
$\pwset{X\times X}$, namely the partial order given by inclusion of binary relations, the multiplication given by
composition, $R\then S=S\circ R$, the
unit $e$ given by the identity relation $\Delta_X$, and the involution given by reversal,
\[R^{*}=\{(y,x)\in X\times X\st xRy\}\;,\]
this mapping yields an isomorphism of quantales from the quantale of sup-preserving mappings from $\pwset X$ to itself to the quantale $\pwset{X\times X}$ of binary relations on the set $X$, leading to the following:

\begin{definition}\label{def:relrep}\label{def:relational.representation}
    Let $X$ be a set. By the \emph{relational quantale} determined by $X$ will be meant the quantale \[\End (X)\] of
    sup-preserving endomorphisms of the orthocomplemented sup-lattice  $\pwset X$, which we will identify henceforth with the quantale $\pwset{X\times X}$ of binary relations on the set $X$. 

    By a \emph{relational representation} of a quantale $Q$ on a set $X$ will be meant a
    homomorphism
    \[r:Q\rightarrow\Rel(X)\]
    from $Q$ to the relational quantale $\Rel(X)$ determined by $X$.

      By a \emph{state} of the representation will be meant an element $x\in X$. For any $a\in Q$, and any states $x,y\in X$, we shall usually write
    \[x\goto a y\]
    to denote that $(x,y)\in r(a)$.

    A representation will be said to be \emph{faithful} provided that it is a homomorphism that is injective.
\end{definition}

Let $Q$ be a quantale and $X$ a set. Then it may be remarked that a mapping
$r:Q\rightarrow\Rel(X)$ is
a relational representation exactly if the following
conditions hold for all $x,y\in X$,
$a,b\in Q$, and $\Phi\subseteq Q$:
\[
\begin{array}{rcl}
    x\goto e x\;,\\
    x\goto{a\then b}y & \iff & x\goto a z\goto b y
    \textrm{ for some }z\in X\;,\\
    x\goto{\V\Phi}y & \iff & x\goto a y\textrm{ for some }a\in\Phi\;,\\
    x\goto{a^{*}}y & \iff & y\goto a x\;.
\end{array}
\]
Furthermore, the representation $r$ is strictly unital if, and only if, the first
condition is strengthened to
\[x\goto e y \iff x=y\]
    for all $x,y\in X$.
From these conditions, and from the fact that $1$ is always
self-adjoint and idempotent, it follows that $r(1)\subseteq
X\times X$ is an equivalence relation
on $X$, and that furthermore $x\goto 1 y$ holds if, and only if,
$x\goto a y$ holds for some $a\in Q$, since $1=\V Q$, leading to the following:

\begin{definition}\label{def:conndet}\label{def:connected}
    Let $Q$ be a quantale, $X$ a set, and let
    $r:Q\rightarrow\Rel(X)$ be a relational representation.
    Then the states $x,y\in X$ will be said to be \emph{connected} in the
    representation provided that \[x\goto 1 y\;.\] The representation itself will be said to be 
    \emph{connected} if \[r(1)=X\times X\;,\] and \emph{deterministic} provided that for each pair of
    connected states $x,y\in X$ there is an element $a\in Q$ such that \[\{x\} a=\{y\}\;.\]
\end{definition}

    It may be remarked that:

\begin{proposition}\label{prop:det.implies.unital}
    Any relational representation \[r:Q\rightarrow\Rel(X)\] that is deterministic is necessarily strictly unital.
\end{proposition}

\begin{proof}
    Let $x$ and $y$ be states of a deterministic representation such
    that $x\goto e y$, and let $a$ be such that $\{y\} a=\{x\}$.
    Then $y\goto{a\then e}y$, hence $y\goto a y$, and thus $x=y$. \qed
\end{proof}

Moreover, it may be noted that for the particular case of relational representations one has the following characterisations of familiar concepts:

\begin{proposition}\label{prop:algebraic.irreducible}
    Let $Q$ be a quantale, $X$ a set, and let
    $r:Q\rightarrow\Rel(X)$ be a relational representation on the set $X$.
    \begin{enumerate}
    \item The representation is irreducible if, and only if,
    it is connected.

    \item The representation is algebraically
    irreducible if, and only if, it is both irreducible and deterministic.
    \item Let
    $r':Q\rightarrow\Rel(X')$ be another relational 
    representation. Then $r$ and $r'$ are equivalent if, and only
    if, there exists 
    a bijective mapping $f:X\rightarrow X'$ such that
    \[x\goto a y\iff f(x)\goto a f(y)\] for all
    for all $x,y\in X$ and $a\in Q$.
    \end{enumerate}
\end{proposition}

\begin{proof}
    1. Connectedness is obviously equivalent to irreducibility, since 
    it states that $r(1)$ equals $X\times X$, that is, 
    $r$ is strong.\\
    2. Equally obviously.\\
    3. If such a bijection exists then its sup-preserving extension
    $f^\sharp:\pwset X\rightarrow\pwset {X'}$ is an isomorphism of 
    sup-lattices, and we obtain, for all $a\in Q$ and $Y\in\pwset X$,
    \[f^\sharp(Ya)=f^\sharp(\bigcup_{y\in Y}\{x\in X\st y\goto a x\}=\]
    \[=\bigcup_{y\in Y}\{f(x)\st f(y)\goto a f(x)\}=
    \bigcup_{w\in f^\sharp(Y)}\{z\in X'\st w\goto a z\}=f^\sharp(Y)a\;;\]
    that is, $f^\sharp$ is an isomorphism of right $Q$-modules. 
    Finally, orthogonality in $\pwset X$ is defined by $Y\perp Z\iff Y\cap 
    Z=\emptyset$, and thus it is preserved and reflected by $f^\sharp$.
    \qed
\end{proof}

    With these preliminaries, we may now consider the decomposition of any relational representation into its irreducible components, noting first the following:

\begin{lemma}\label{lem:canonical.decomposition}
    Let $(X_i)_{i\in I}$ be a family of sets indexed by a set $I$. Then
    the direct product \[Q=\prod_{i\in I}\Rel(X_{i})\]
    has a canonical representation on the disjoint union
    $X=\coprod_{i\in I}X_{i}$, which is defined by
    inclusion and is faithful.
\end{lemma}

\begin{proof}
    Assume for simplicity that the sets $X_{i}$ are pairwise disjoint
    (if not, turn the representations into equivalent ones
    by labelling the states of each $X_{i}$ with $i$),
    let $X=\bigcup_{i\in I} X_{i}$,
    and define a map $\kappa: Q\rightarrow\Rel(X)$ by
    $\kappa(a)=\bigcup_{i\in I} R_{i}$, for all $a=(R_{i})_{i\in I}\in
    Q$. Then $\kappa$ clearly is strictly unital, it preserves joins,
    and it preserves
    involution. It also preserves multiplication, which is most
    easily seen by representing binary relations as matrices with
    entries in the two element chain (\ie, matrices of zeros and
    ones with addition replaced
    by the operation of taking binary joins) and $\kappa$ as the map
    that sends each family of
    matrices as blocks into a large matrix indexed by $X\times X$;
    then multiplication being preserved just means that the multiplication
    of two such block matrices is computed blockwise. \qed
\end{proof}

\begin{definition}\label{def:sum.of.reps}
    Let $Q$ be a quantale, and let
    $(r_{i}:Q\rightarrow\Rel({X_{i}}))_{i\in I}$ be a
    family of relational representations indexed by a set $I$. By the
    \emph{sum} of the representations, denoted by \[r=\coprod_{i\in
    I}r_{i}\;,\] will be meant the relational representation of $Q$ on the disjoint union
    $\coprod_{i\in I}X_{i}$ defined by
    \[r(a)=\kappa((r_{i}(a))_{i\in I})\;,\]
    where $\kappa$ is the canonical representation of
    $\prod_{i\in I}\Rel({X_{i}})$ on $\coprod_{i\in I}X_{i}$.
\end{definition}

These definitions now lead to the following evident conclusion:

\begin{theorem}\label{thm:irredreps}
    Let $Q$ be a quantale with a relational representation 
    \[r:Q\rightarrow\Rel(X)\] on a set $X$.
    Then $r:Q\rightarrow\Rel(X)$ admits a canonical decomposition into a sum of irreducible
    representations, one on each connected component of $X$. Furthermore
    this decomposition is unique up to equivalence of representations,
    and the irreducible representations are all algebraically
    irreducible if, and only if, the relational representation $r:Q\rightarrow\Rel(X)$ is deterministic.
\end{theorem}

    It may be noted in passing that this situation is evidently a particular aspect of relational representations, which simplifies considerably our discussion from this point.

\section{The relational representations of $\penq$}\label{sec:5}

In this section we begin to study the points of the quantale $\penq$. Contrary to the situation for quantales of the form
$\Q A$, where we are concerned with points that are representations on sup-lattices of projections of Hilbert spaces, we shall see that for $\penq$ the ``natural" points are relational representations. Indeed, these are the only points we shall study, in particular concluding that they correspond bijectively with Penrose tilings. No other points of the quantale $\penq$ are known.

It should be emphasised that, because of its particular properties, at one time both very different from quantales of the form
$\Q A$ but closely related to them through the C*-algebra $\pena$ introduced by Connes, our interest is almost as much in the light that the quantale $\penq$ can throw on the concept of point as on the information that points can contribute to the study of the theory $\pent$ of Penrose tilings. For this reason, we begin by examining more closely the weaker concept of an irreducible representation
of the quantale $\penq$ on the powerset of a set $X$, building concretely on the results of the previous section.

The first step will be to show that for any relational representation \[r:\penq\rightarrow\Rel(X)\]
of the quantale $\penq$ there is a natural way in which a Penrose sequence $\seq(x)$ can be assigned to each state $x$ of the representation (Definition~\ref{def:tiling.sequence}). In order to do this we need the following preliminary lemma, where in order to simplify notation we omit angular brackets in 
expressions like $x\goto{\wprop n X}y$, writing instead
$x{\wprop n X}y$ to mean that $(x,y)\in r(\wprop n 
X)$ for a representation $r$ of $\penq$. For instance, this 
gives us the following convenient equivalence:
\[x{\wprop n X}y\iff y{\rprop n X}x\;.\] These simplifications will be used without further comment through the rest of the paper.

\begin{lemma}\label{lem:tiling.sequence}
    Let $r:\penq\rightarrow\Rel(X)$ be a relational 
    representation of the quantale $\penq$ on a set $X$. Then for any $x\in X$ and all $n\in\mathbb N$, the following conditions are equivalent:
    \begin{enumerate}
	\item $x{\wprop n L} x$;
	\item not $x{\wprop n S} x$;
	\item $y{\wprop n L} x$ for some $y\in X$;
	\item $y{\wprop n S} x$ for no $y\in X$.
    \end{enumerate}
    Furthermore, if $x{\wprop n S} x$ then
    $x{\wprop {n+1} L} x$.
\end{lemma}

\begin{proof}
    Since for any $x\in X$ we always have $x\goto e x$, the axiom $D1_{n}$
    of Definition~\ref{def:pent} tells us that either $x{\wprop n L} x$ or
    $x{\wprop n S} x$. If we had both $x{\wprop n L}x$ and
    $x{\wprop n 
    S}x$ then we would also have $x{\rprop n S}x$ and
    $x{\langle\bprop n L\vert\bprop n S\rangle} x$, which is 
    impossible due to axiom $C1_{n}$. Hence, the first two conditions are 
    equivalent. The first condition implies the third trivially, and 
    the third implies the first because $y{\wprop n L} 
    x{\wprop n S} x$ implies $y{\langle\bprop n L\vert\bprop 
    n S\rangle}x$, which is impossible, again due to $C1_{n}$,
    and thus if $y{\wprop n 
    L} x$ we must have $x{\wprop n L}x$. In a similar way we prove 
    that the second and the fourth conditions are equivalent to each 
    other. Finally, $x{\wprop n S} x$ entails, by axiom $C2_{n}$,
    the condition $x{\wprop{n+1} L}x$. \qed
\end{proof}

\begin{definition}\label{def:tiling.sequence}
    Let $X$ be a set, let $r:\penq\rightarrow\Rel(X)$ be a 
    relational representation, and let $x\in X$. Then by the 
    \emph{tiling sequence} of $x$ will be meant the Penrose sequence $\seq(x)$
    defined by
    \[\seq(x)_{n}=\left\{\begin{array}{ll}
    0 & \textrm{if }x{\wprop n L} x\\
    1 & \textrm{if }x{\wprop n S} x
    \end{array}\right.\]
for each $n\in\mathbb N$.
\end{definition}

Now we show, in the following two lemmas, that the  transitions induced by primitive propositions,
$x\wprop n X y$, are determined entirely by the sequences $\seq(x)$ and $\seq(y)$, provided that the states $x$ and $y$ are connected (\ie, that $x\goto 1 y$).

\begin{lemma}\label{lem:connectivity.by.generators}
    Let $X$ be a set, let $r:\penq\rightarrow\Rel(X)$ be a 
    relational representation, and let $x,y\in X$. The following 
    conditions are equivalent:
    \begin{enumerate}
	\item $x\goto 1 y$;
	\item $x{\wprop n X} y$ for some generator $\wprop n X$.
    \end{enumerate}
    Furthermore, if $x{\wprop n X} y$ then 
    $\seq(x)_{k}=\seq(y)_{k}$ for all $k>n$. (Hence, 
    $\seq(x)\sim\seq(y)$.)
\end{lemma}

\begin{proof}
    The second condition implies the first because $\wprop m Z\le 1$ for any $Z$ and $m$. 
    So let us assume that $x\goto 1 y$.
    This condition means that for some string
    $x_{1}\ldots x_{r}$, where each $x_{i}$ is a quantale generator
    $\wprop m Z$ or $\rprop m Z$, we have $x\goto a y$ for 
    $a=x_{1}\then\ldots\then x_{r}$. Define the 
    \emph{degree} $\deg(x_{1}\ldots x_{r})$ of each such string to be
    $\max\{\deg(x_{1}),\ldots,\deg(x_{r})\}$, where for the generators
    we define $\deg(\wprop m Z)=\deg(\rprop m Z)=m$.
    Let then $x\goto a y$ for $a$ as above, and let $n=\deg(a)+1$. We must
    have $y{\wprop n X} y$ for $X$ equal either to $S$ or $L$, 
    and thus $x\goto{a\then\wprop n X} y$. But from the expansion 
    axioms of Definition~\ref{def:pent} it follows that $a\then\wprop n 
    X\le\wprop n X$ 
    in $\penq$, and thus $x\wprop n X y$, thus proving $1\Rightarrow 2$.
    Now let $k>n$. Then $y\rprop k Y y$ for some $Y$,
    whence $x\langle\bprop n X\vert\bprop k Y\rangle y$. By repeated 
    application of one of the expansion axioms we obtain
    $\langle\bprop n X\vert\bprop k Y\rangle\le\rprop k Y$ in $\penq$, and 
    thus $x\rprop k Yy$. This implies that $x\wprop k Y x$, 
    for otherwise we would have $x\langle \bprop k Z\vert\bprop 
    k Y\rangle y$ with $Z\neq Y$, which is impossible because
    $\langle \bprop k Z\vert\bprop k Y\rangle=0$ if $Z\neq Y$, due to one 
    of the consistency axioms. Hence, we conclude that $\seq(x)_{k}=\seq(y)_{k}$. \qed
\end{proof}

\begin{lemma}\label{lem:characterizing.connectivity}
    Let $X$ be a set, let $r:\penq\rightarrow\Rel(X)$ be a
    relational representation, and let $x,y\in X$ be such that $x\goto 1 y$.
    Then, we have $x{\wprop n L} y$ (resp.\ $x{\wprop n S} 
    y$) if, and only if, we have both $\seq(y)_{n}=0$ (resp.\ $\seq(y)_{n}=1$)
    and $\seq(y)_{k}=\seq(x)_{k}$ for all $k>n$.
\end{lemma}

\begin{proof}
    The forward implication is an immediate consequence of the 
    previous lemma. For the reverse implication, assume that
    we have both $\seq(y)_{n}=0$ (resp.\ $\seq(y)_{n}=1$)
    and $\seq(y)_{k}=\seq(x)_{k}$ for all $k>n$. Since we are 
    assuming that $x\goto 1 y$, we have
    $x{\wprop m X} y$ for some generator $\wprop m X$. If $m=n$ 
    then $X$ must be $L$ (resp.\ $S$) due to one of the consistency 
    axioms, and we are done. If $m<n$ we conclude that 
    $x{\wprop n L} y$ (resp.\ $x{\wprop n S}y$) by applying
    one of the expansion axioms (as in the proof of 
    the previous lemma). Now assume that $m>n$. Then 
    $\seq(x)_{m}=\seq(y)_{m}$, and thus $x{\rprop m X} x$.
    Hence, we obtain $x{\vert\bprop m X\rangle\langle \bprop m X\vert} y$,
    and by applying the second decidability axiom we conclude that
    $x{\wprop{m-1} L} y$ or $x{\wprop{m-1} S} y$ (the 
    choice of which is determined by the value of $\seq(y)_{m-1}$). 
    The proof follows by induction on $m-n$. \qed
\end{proof}

The previous two lemmas describe completely the transitions
$x\goto a y$ of $\penq$ between any pair of connected states $x$ and $y$ of the set $X$. 
The following lemma improves this by showing that there must be ``as many 
states as possible'', \ie, each irreducible component must be saturated with respect to equivalence of Penrose sequences.

\begin{lemma}\label{lem:enough.states}
    Let $X$ be a set, let $r:\penq\rightarrow\Rel(X)$ be 
    a nonempty 
    relational representation, and let $x\in X$.
    Then for any Penrose sequence
    $s$ such that 
    $s\sim\seq(x)$, there must be a state $y\in X$ such that 
    $s=\seq(y)$.
\end{lemma}

\begin{proof}
    Let $s_{n}=\seq(x)_{n}$ for all $n>m\in\mathbb N$, and consider 
    the product of quantale generators
    $\langle\bprop m{X_{m}};\ldots;\bprop 0{X_{0}}\vert
    \bprop 0{X_{0}};\ldots;\bprop m{X_{m}}\rangle$, where for 
    each $i\in\{0,\ldots,m\}$ we have $X_i=L$ if, and only if, $s_{i}=0$. 
    Without loss of generality we may assume that $X_{m}=L$ (if not, 
    replace $m$ by $m+1$),
    and the completeness axiom tells us that
    $(x,x)\in r(\langle\bprop m{X_{m}};\ldots;\bprop 0{X_{0}}\vert
    \bprop 0{X_{0}};\ldots;\bprop m{X_{m}}\rangle)$. Hence, 
    there exists a state $y\in X$ such that $(x,y)\in
    r(\langle\bprop m{X_{m}};\ldots;\bprop 0{X_{0}}\vert)$, and from 
    Lemma~\ref{lem:characterizing.connectivity} it is immediate that $s_{i}=\seq(y)_{i}$ 
    for all $i\in\{0,\ldots,m\}$, and thus for all $i\in\mathbb N$. 
    \qed
\end{proof}

Taking these results together, we have established the following result:

\begin{theorem}\label{thm:seq}
Any relational representation $r:\penq\rightarrow\Rel(X)$ of the quantale $\penq$ on a set $X$ determines a mapping
\[\seq:X\rightarrow K\]
from the set $X$ to the set $K$ of Penrose sequences by assigning to each $x\in X$ the sequence defined by
    \[\seq(x)_{n}=\left\{\begin{array}{ll}
    0 & \textrm{if }x{\wprop n L} x\\
    1 & \textrm{if }x{\wprop n S} x
    \end{array}\right.\]
for each $n\in\mathbb N$, in such a way that each connected component of $X$ is mapped surjectively onto an equivalence class of Penrose sequences.
\end{theorem}

Applying this to the case of an irreducible representation, characterised by the set $X$ having a single connected component, we have the following observation:

\begin{corollary}\label{cor:enough.states}
    For any irreducible representation of $\penq$ on a set $X$, one has that:
\begin{enumerate}
\item the subset $\{\seq(x)\in K\st x\in X\}$ is an equivalence class of Penrose sequences, and
\item the action of the quantale $\penq$ on the set $X$ is such that, for any pair of states $x,y\in X$, we have
$x{\wprop n L} y$ (resp.\ $x{\wprop n S} 
    y$) if, and only if, both $\seq(y)_{n}=0$ (resp.\ $\seq(y)_{n}=1$)
    and $\seq(y)_{k}=\seq(x)_{k}$ for all $k>n$.
\end{enumerate}
\end{corollary}

To obtain the converse, giving a complete characterisation of the relational representations of $\penq$ in terms of Penrose sequences, we first make the following:

\begin{definition}
Let $X$ be a set and $(X_i)_{i\in I}$ a family of subsets that partition the set $X$. Suppose given a mapping
\[\sigma:X\rightarrow K\]
such that the image of each subset $X_i$ of the partition is an equivalence class of Penrose sequences. By the \emph{interpretation induced by $\sigma$} of the theory $\pent$ in the quantale $\Rel(X)$ is meant that obtained by assigning to each of the primitive propositions the relations defined by:
    \begin{eqnarray*}
	\wprop n L_{\sigma} &=& \{(x,y)\in \Xi\st \sigma(y)_{n}=0\textrm{ and }\sigma(x)_{k}=
	\sigma(y)_{k}\textrm{ for all }k>n\}\;,\\
	\wprop n S_{\sigma} &=& \{(x,y)\in \Xi\st \sigma(y)_{n}=1\textrm{ and }\sigma(x)_{k}=
	\sigma(y)_{k}\textrm{ for all }k>n\}\;,
    \end{eqnarray*}
where $\Xi=\bigcup_{i\in I} X_i\times X_i$ is the equivalence relation determined by the partition.
     \end{definition}

Showing that this interpretation validates the axioms of the theory $\pent$ in the quantale $\Rel(X)$ to obtain a relational representation of the quantale $\penq$ on the set $X$ then yields the following converse to the result already obtained:

\begin{theorem}\label{thm:characterize.irred.reps}
Let $X$ be a set and $(X_i)_{i\in I}$ a family of subsets that partition the set $X$. Let \[\sigma:X\rightarrow K\] be a mapping for which
the image of each subset $X_i$ of the partition is an equivalence class of Penrose sequences.
    Then the interpretation of the primitive propositions induced by $\sigma$ in $\Rel(X)$ extends in a unique way to an irreducible relational representation
    $r:\penq\rightarrow\Rel(X)$
for which \[\seq(x)=\sigma(x)\]
    for each $x\in X$. Moreover, any relational representation of $\penq$ arises uniquely in this way.
\end{theorem}

\begin{proof}
    The primitive propositions may be considered to be the generators of $\penq$, so to obtain the required homomorphism
    we have to verify that the 
    defining relations of Definition~\ref{def:pent} are respected, \ie, that each of the axioms of the theory $\pent$ is validated by the given interpretation.

    ($C1_{n}$) Proving that this relation is respected means showing that
    \[\wprop n L_\sigma\then \rprop n S_\sigma\subseteq\emptyset\;.\]
    Let then  $x\wprop n L_\sigma y\rprop n S_\sigma z$. We have both
    $x\wprop n L_\sigma y$ and $z\wprop n S_\sigma y$, and thus both
    $\sigma(y)_{n}=0$ and 
    $\sigma(y)_{n}=1$, a contradiction. Hence, $x\langle\bprop n L\vert \bprop n S\rangle_\sigma z$ for no two states $x$ and $z$.

    ($C2_{n}$) Let $x\wprop n S_\sigma y$ (hence, $x$ and $y$ are both in the same component of the partition). Then we have $\sigma(y)_{n}=1$ and 
    $\sigma(x)_{m}=\sigma(y)_{m}$ for all $m>n$. Then also $\sigma(y)_{n+1}=0$, because 
    $\sigma(y)$ is a Penrose sequence, and $\sigma(x)_{m}=\sigma(y)_{m}$ for all $m> 
    n+1$, whence $x\wprop {n+1}L_\sigma y$ (because $x$ and $y$ are in the same component of the partition).
    Hence, we conclude that
    $\rprop n S_\sigma\subseteq \wprop{n+1}L_\sigma$.

    ($D1_{n}$) Let 
    $x\in X$. Then $\sigma(x)_{n}=0$ or $\sigma(x)_{n}=1$, and thus
    $x\wprop n L_\sigma x$ or $x\wprop n S_\sigma x$. This is 
    equivalent to $(x,x)\in \wprop n L_\sigma\cup \wprop n S_\sigma$, showing 
    that the relation is respected.

    ($D2_{n}$) Let $(x,y)\in\rprop {n+1}X_\sigma\then\wprop {n+1}X_\sigma$. We have
    $\sigma(x)_{n+1}=\sigma(y)_{n+1}=0$ if $X=L$, and
    $\sigma(x)_{n+1}=\sigma(y)_{n+1}=1$ if $X=S$. In either case, we have 
    $\sigma(x)_{n+1}=\sigma(y)_{n+1}$, and thus $\sigma(x)_{m}=\sigma(y)_{m}$ for all
    $m>n$. Then we either have $x\wprop n S_\sigma y$, if 
    $\sigma(y)_{n}=1$, or $x\wprop n L_\sigma y$, if 
    $\sigma(y)_{n}=0$. Hence, we conclude that
    $\rprop {n+1}X_\sigma\then \wprop {n+1}X_\sigma\subseteq \wprop n L_\sigma\cup \wprop n S_\sigma$.

    ($E1_{n}$) Let us assume for example that $X=L$. If
    $(x,y)\in \wprop n Y_\sigma\then\wprop {n+1}L_\sigma$ then we clearly have 
    both $\sigma(y)_{n+1}=0$ and $\sigma(y)_{m}=\sigma(x)_{m}$ for all $m>n+1$, and 
    thus $x\wprop{n+1}L_\sigma y$. If $X=S$ the conclusion is similar.

    ($E2_{n}$)--($E4_{n}$) These are handled in a similar way to the previous
    one.

    ($I1_{n}$) This just states that $x\rprop n X_\sigma y$
    implies $y\wprop n X_\sigma x$ .

    ($I2_{n}$) Similarly.

    ($C'_{t}$) Let $X_i$ be an arbitrary subset in the partition and let $x\in X_i$. Consider a string
    $\bprop n{X_{n}};\ldots;\bprop 0{X_{0}}$, with each $X_{k}$ equal 
    to $S$ or $L$. If this string is 
    admissible in the sense of Definition~\ref{def:pent}, and $X_{n}=L$,
    then the sequence $t:\mathbb N\rightarrow\{0,1\}$ defined by
    \[t_{k}=\left\{\begin{array}{ll}
    \sigma(x)_{k} &\textrm{if }k>n\\
    0 &\textrm{if }k=n\\
    0 &\textrm{if }k<n\textrm{ and }X_{k}=L\\
    1 &\textrm{if }k<n\textrm{ and }X_{k}=S
    \end{array}\right.\]
    is a Penrose sequence, which furthermore is equivalent to $\sigma(x)$. 
    Hence, there must exist $y\in X_i$ such that $\sigma(y)=t$ because the 
    image of $X_i$ under $\sigma$ is a whole equivalence class of Penrose sequences, and
    we obtain the following conditions:
    \[\begin{array}{l}
	x\wprop n{X_{n}}_\sigma y\;,\\
	y\wprop k{X_{k}}_\sigma y\textrm{ for all }k<n\;.
    \end{array}\]
    As a consequence we conclude that
    $x\langle\bprop n{X_{n}};\ldots;\bprop 0{X_{0}}\vert_\sigma y$, 
    and finally that
    $x\langle\bprop n{X_{n}};\ldots;\bprop 0{X_{0}}\vert
    \bprop 0{X_{0}};\ldots;\bprop n{X_{n}}\rangle_\sigma x$, thus showing 
    that
    \[\Delta_X\subseteq \langle\bprop n{X_{n}};\ldots;\bprop 0{X_{0}}\vert
    \bprop 0{X_{0}};\ldots;\bprop n{X_{n}}\rangle_\sigma\;. \]
    In consequence, the interpretation determined by the mapping $\sigma:X\rightarrow K$ over the given partition
    $(X_i)_{i\in I}$ of $X$ determines a relational representation of $\penq$ on the set $X$ having the required property. The converse, that any such representation arises in this way from a unique mapping $\sigma$, is of course a consequence of
    Theorem~\ref{thm:seq}. \qed
    \end{proof}

It may be remarked that of course any mapping $\sigma:X\rightarrow K$ of which the image is a disjoint union of equivalence classes of Penrose sequences gives rise to a canonical partition of the set $X$, namely by taking the inverse images of the equivalence classes of Penrose sequences. However, not every relational representation of the quantale $\penq$ arises in this way, but only those whose canonical decomposition into irreducible representations admits at most one instance of each equivalence class of irreducible representations. It is to allow multiplicity of irreducible components that the characterisation of relational representations takes this more general form.

In the particular case of irreducible representations, we have therefore the following characterisation:

\begin{corollary}
Any irreducible representation \[r:\penq\rightarrow\Rel(X)\]
of the quantale $\penq$ on a set $X$ determines, and is uniquely determined by, a mapping
\[\sigma:X\rightarrow K\]
of which the image is an equivalence class of Penrose sequences in such a way that $\seq(x)=\sigma(x)$ for each $x\in X$.
\end{corollary}

\section{The relational points of $\penq$}\label{sec:6}

In this section we finally describe the relational points of the quantale $\penq$, that is to say, the algebraically irreducible representations of $\penq$ on the powerset lattices of sets. More than that, we see the way in which these points emerge from the motivating description of the axioms for the noncommutative theory $\pent$ of Penrose tilings that gives rise to the quantale $\penq$, along with the Penrose sequences with which we have been working.

To begin with, we examine the characterisation of relational representations obtained in the preceding section, in the case of the canonical interpretation given in terms of translations that was considered informally in introducing the theory $\pent$ in section~\ref{sec:3}, showing that this yields a natural description of the irreducible representations to which it gives rise.

\begin{example}\label{exm:geometrical.rep}
The canonical interpretation leads naturally to a relational representation of the quantale $\penq$ on the set $\Sigma$ of all Penrose tilings for which the origin of the plane does not lie on a vertex or an edge of any tile, by interpreting the primitive propositions of the theory in the way that we shall now describe. Denoting by $T$ and $T'$ the tiles (of level $0$) of $\mathbb T$ and $\mathbb T'$, respectively, for which $0\in T$ and $0\in T'$, we write
$^n T$ and $^n T'$ for their $n$th successor tiles $s^n(T)$ and $s^n(T')$, respectively: that is to say, the unique tiles of level $n$ that contain the origin in each tiling. Then, each primitive proposition of the theory $\pent$ is considered to describe transitions between states $\mathbb T,\mathbb T'\in\Sigma$, defined by:
\begin{itemize}
\item $\mathbb T{\wprop n L}\mathbb T'$ if, and only if, $\mathbb T'=\mathbb
T+z$ for some $z\in\mathbb C$ such that $0\in s({}^n T+z)$ and ${}^n T'$ 
is large (at level $n$);
\item $\mathbb T{\wprop n S}\mathbb T'$ if, and only if, $\mathbb T'=\mathbb
T+z$ for some $z\in\mathbb C$ such that $0\in s({}^n T+z)$ and ${}^n T'$ is small;
\item $\mathbb T{\rprop n L}\mathbb T'$ if, and only if, ${}^n T$ is large
and $\mathbb T'=\mathbb T+z$ for some $z\in\mathbb C$ such that $0\in
s({}^n T+z)$;
\item $\mathbb T{\rprop n S}\mathbb T'$ if, and only if, ${}^n T$ is small and
$\mathbb T'=\mathbb T+z$ for some $z\in\mathbb C$ such that $0\in
s({}^n T+z)$.
\end{itemize}
The axioms of the theory $\pent$ of Penrose tilings are validated by this canonical interpretation, yielding in consequence a relational representation of the quantale $\penq$ on the set $\Sigma$
whose restriction to primitive propositions is that described above: we shall refer to this as the \emph{geometrical representation},
\[g:\penq\rightarrow\Rel(\Sigma)\;,\]
of the quantale $\penq$.

In terms of the characterisation of relational representations of the preceding section, this representation is classified by the mapping
\[\gamma:\Sigma\rightarrow K\]
obtained by assigning to each tiling the Penrose sequence generated by the origin of the plane, together with the canonical partition $(\Sigma_i)_{i\in I}$ of the set $\Sigma$ induced by the mapping, namely by taking the inverse image of each equivalence class of Penrose sequences.
\end{example}

It is of course the case in the above example that each subset $\Sigma_i$ of the canonical partition contains tilings of which the equivalence is evident in different ways: some directly from translations describable in terms of the experiments represented by the primitive propositions, some dependent additionally on other operations, such as rotations of the plane, which have not been coded directly into the experiments concerned. In other words, the geometrical representation $g:\penq\rightarrow\Rel(\Sigma)$ to which the canonical interpretation gives rise will decompose uniquely into irreducible representations $g_i:\penq\rightarrow\Rel(\Sigma_i)$ determined by sets, consisting of all tilings of a particular equivalence class, that contain tilings that are not related deterministically by the experiments considered. In terms of states, this means that we have an equivalence of states that arises from symmetries that have not been addressed. In terms of the representation, this means that while it decomposes into irreducible representations, these irreducible representations are not actually algebraically irreducible, a situation that is characterised by the following:

\begin{theorem}\label{thm:algirred.iff.irredun}
    Let $X$ be a set, and let $r:\penq\rightarrow\Rel(X)$ be
    an irreducible 
    relational representation on the set $X$. Then $r$ is algebraically irreducible
    if, and only if, the canonical mapping
    $\seq:X\rightarrow K$ which it determines is injective.
\end{theorem}

\begin{proof}
    Assume that $r$ is algebraically irreducible, and
    let $x,y\in X$ be such that $\seq(x)=\seq(y)$. Let also $a\in\penq$ 
    be a product of generators $a=\alpha_{1}\then\ldots\then 
    \alpha_{n}$. We prove that if $x\goto a y$ then $x\goto a x$, by 
    induction on $n$.
    If $n=0$ (\ie, $a=e$), then this is trivial because $x\goto e y$ 
    implies $x=y$, since any algebraically irreducible representation 
    is strictly unital. Let then $n>0$, and $x\goto a y$. We have
    $a=b\then \alpha_{n}$, and 
    thus $x\goto b z\goto{\alpha_{n}} y$ for 
    some $z\in X$. From Theorem~\ref{thm:characterize.irred.reps} 
    we conclude that $z\goto{\alpha_{n}} x$ because $z\goto{\alpha_n}y$ and $\seq(x)=\seq(y)$. Hence, 
    $x\goto a x$. Let now $a\in\penq$ be such that
    $\{x\} a=\{y\}$; such an $a\in\penq$ exists because $r$ is 
    deterministic. Without loss of generality we may assume that $a$ 
    is a product of generators, and thus from $x\goto a y$ we 
    conclude $x\goto a x$, which implies $x\in\{x\}a=\{y\}$, \ie, $x=y$. 
    Hence, $\seq$ is injective.

    Now we prove that $r$ is algebraically irreducible
    if $\seq$ is injective. Let then $x$ and $y$ be two states of the representation, and let $s=\seq(x)$ and $t=\seq(y)$. Since the representation is irreducible, the Penrose sequences $s$ and $t$ must be equivalent. Let $n\in\mathbb N$ be such that $s_{k}=t_{k}$ for all 
    $k\ge n$. Without loss of generality assume that $s_{n}=t_{n}=0$.
    Now let $a=\langle \bprop n{X_{n}};\ldots;\bprop 0{X_{0}}\vert$, 
    where for each $k\in\{0,\ldots,n\}$ we have
    $X_{k}=L$ if, and only if, $t_{k}=0$. Then clearly we have $x\goto 
    a y$, and furthermore if $x\goto a z$ we must have 
    $\seq(z)=t=\seq(y)$. Hence, $z=y$ because $\seq$ is injective, and therefore we have
    $\{x\}a=\{y\}$, thus showing that the representation $r$ is algebraically 
    irreducible. \qed
\end{proof}

Of course, an immediate consequence is the following characterisation of the relational representations of the quantale
$\penq$ that are algebraically irreducible, hence of the relational points of the quantale of Penrose tilings:

\begin{corollary}
A relational representation $r:\penq\rightarrow\Rel(X)$ is algebraically irreducible if, and only if,
it is equivalent to the relational representation
\[\rho_x:\penq\rightarrow\End(x)\]
canonically determined by an equivalence class $x\in K/{\sim}$ of Penrose sequences, by which we mean the unique representation on $x$ for which the canonical mapping $\seq:x\rightarrow K$ is the inclusion of $x$ into $K$. In particular, the relational points of the quantale
$\penq$ are in canonical bijective correspondence with the equivalence classes of Penrose sequences.
\end{corollary}

The effect of moving to the consideration of Penrose sequences in studying Penrose tilings may therefore be seen to be precisely that of taking us away from the redundancies intrinsically involved in considering irreducible representations of the quantale $\penq$ to the consideration of algebraically irreducible representations, in other words to the geometric content of the quantale $\penq$ of Penrose tilings described by its points.

In terms of the classification of relational representations of the quantale $\penq$, the sum of the algebraically irreducible representations of $\penq$ determined by its points is exactly that classified by the identity mapping
\[\mathrm{id}_K:K\rightarrow K\]
together with the canonical partition of the set $K$ into equivalence classes of Penrose sequences, yielding the following:

\begin{definition}\label{def:Penrose.rep}
    By the \emph{Cantor representation} of $\penq$ will be meant the relational 
    representation whose set of states is $K$, together with the interpretation of primitive propositions defined as follows:
\begin{itemize}
\item $s\wprop n L t$ if, and only if, both $t_{n}=0$ and
    $s_{m}=t_{m}$ for all $m>n$;
\item $s\wprop n S t$ if, and only if, both $t_{n}=1$ and
    $s_{m}=t_{m}$ for all $m>n$.
\end{itemize}
\end{definition}

The Cantor representation is canonical in the following sense:

\begin{theorem}
    Let $X$ be a set, and let $r:\penq\rightarrow\Rel(X)$ be
    a relational representation. Then the map $\seq:X\rightarrow K$ 
    extends canonically to a homomorphism of right $\penq$-modules \[\Seq:\pwset 
    X\rightarrow\pwset K\;.\]
\end{theorem}

\begin{proof}
    It suffices that we show, for each generator $a=\wprop n X$ or $a=\rprop n 
    X$, and for each $x\in X$, that $\Seq(\{x\}a)=\Seq(\{x\})a$, 
    that is, $\{\seq(y)\in K\st x\goto a y\}=\{s\in K\st \seq(x)\goto a 
    s\}$,
    which means proving for all $x\in X$ and all $s\in K$ that
    we have $\seq(x)\goto a s$ in the Cantor representation if, and only if, 
    $x\goto a y$ for some $y\in X$ such that $\seq(y)=s$. From 
    Theorem~\ref{thm:characterize.irred.reps} one of the implications 
    is obvious: if $x\goto a y$ then $\seq(x)\goto a\seq(y)$. 
    Furthermore, again by Theorem~\ref{thm:characterize.irred.reps}, 
    each irreducible component maps surjectively to an 
    irreducible component of the Cantor representation, and thus if 
    $\seq(x)\goto a s$ there is $y\in X$ such that $\seq(y)=s$, and 
    furthermore we have $x\goto a y$. \qed
\end{proof}

In particular, the relationship between the geometrical representation and the Cantor representation is that expected from its r\^{o}le in eliminating states that are equivalent for geometric reasons:

\begin{corollary}
    The Cantor representation is a right $\penq$-module
    quotient of the geometrical 
    representation.
\end{corollary}

Of course, the Cantor representation is also canonical
in the sense already described in terms of the relational points of the quantale, but which may equally be stated in terms of its irreducible components in the following way:

\begin{corollary}\label{cor:classif.algirred}
    The Cantor representation is deterministic, and a relational representation
    of $\penq$ is algebraically 
    irreducible if, and only if, it is equivalent to an irreducible 
    component of the Cantor representation.
\end{corollary}

We conclude this section remarking that in view of the canonical r\^{o}le played by Penrose sequences, one may observe that the Cantor representation may be considered to establish an equivalence between the geometrically motivated theory $\pent$ of Penrose tilings and the more algebraically motivated theory defined in the following way:

\begin{definition}\label{def:pens}
By the \emph{theory $\pens$ of Penrose sequences modulo equivalence} will be meant that obtained by taking primitive propositions of the form
\[(s_n=0)\textrm{ and }(s_n=1)\]
for each $n\in\N$, together with the following axioms,
\[
\begin{array}{ll}
(C1_{n}) &(s_n=0)\then (s_n=1)^*\entails
\false\\
(C2_{n}) &(s_n=1)\entails(s_{n+1}=0)\\
(D1_{n}) &\true\entails(s_n=0)\vee(s_n=1)\\
(D2_{n}) &(s_{n+1}=b)^*\then(s_{n+1}=b)\entails
(s_n=0)\vee(s_n=1)\\
(E1_{n}) &(s_n=b')\then(s_{n+1}=b)\entails (s_{n+1}=b)\\
(E2_{n}) &(s_{n+1}=b)\then(s_n=b')\entails (s_{n+1}=b)\\
(E3_{n}) &(s_n=b')^*\then(s_{n+1}=b)\entails(s_{n+1}=b)\\
(E4_{n}) &(s_{n+1}=b)\then(s_n=b')^*\entails(s_{n+1}=b)\\
(C'_{t}) &\true\entails (s_n=b_n)\then\ldots\then(s_0=b_0)\then(s_0=b_0)^*\then\ldots\then(s_n=b_n)^*
\end{array}
\]
for each $b,b'\in\{0,1\}$, and all admissible strings $t=b_0\ldots b_n$ having $b_n=0$, where a string
$b_0\ldots b_n$ of $0$'s and $1$'s is said to be \emph{admissible} provided that it contains no consecutive $1$'s.
\end{definition}

Evidently, this equivalence of theories expresses that the primitive propositions $\wprop n L$ and $\wprop n S$ of the theory $\penq$ of Penrose tilings may be identified with primitive propositions $(s_n=0)$ and $(s_n=1)$, respectively, that describe the computational properties of Penrose sequences. In turn, the quantales determined by these noncommutative theories may each be considered as the noncommutative space $\penq$ of Penrose tilings, of which we have shown that the relational points
are indeed the Penrose tilings of the plane.

\section{Conclusion}\label{sec:7}

In this paper we have obtained the quotient model of the set of Penrose tilings in terms of the representation theory of quantales, taking as a starting point a geometrical description of the tilings by means of a noncommutative theory. Specifically, we have provided a complete classification of the relational points of the quantale $\penq$, showing that they can be identified with the Penrose tilings.
This work falls into the general effort of understanding the r\^{o}le of quantales in providing a generalised notion of space, and it is the first example where \emph{relational} algebraically irreducible representations of a quantale are studied.
It would be therefore interesting to know whether there are other points of $\penq$ besides the relational ones, or to know what other axioms are necessary in the theory $\pent$ in order to rule out nonrelational points. This question is also interesting from a philosophical viewpoint, for it may be argued that the importance of relational representations in this case is directly related to the fact that underlying our axioms is a ``classical" (as opposed to quantum) notion of disjunction (\cf\ section~\ref{sec:3}).

Incidentally to so doing, we have proved some new results about relational representations, in particular obtaining a decomposition theorem for them. Concretely, any relational representation on a set $X$ is partitioned into irreducible components by the connectivity equivalence relation on $X$, and it is this equivalence relation which, in the case of the quantale of Penrose tilings, coincides with the equivalence relation on Penrose sequences that yields as a quotient the set of tilings. Hence, not only the set $K$ of Penrose sequences is derived from the representation theory of $\penq$, but, in addition, and in contrast with the situation for locales, also the equivalence relation on $K$ is derived from $\penq$. This fact makes relational representations interesting in their own right, and meriting further study.

In terms of the theory of quantales, the analysis of relational representations has also provided the opportunity to study, at least in this context, the relationship between irreducible and algebraically irreducible representations of quantales. Explicitly, we have seen that in this context each irreducible representation may be refined to an algebraically irreducible representation by passing to a quotient set of the set of states. The subtleties of this in more general situations will be examined elsewhere.

It has also been remarked that the identification of Penrose tilings with the relational points of the quantale $\penq$ allows us to view the set of tilings as having a quantum aspect, in the sense that the logical assertions about them introduce nondeterministic translational modifications in the tiling being observed, although this aspect is only partial, to the extent that superposition of tiling states by means of linear combination is still absent. In this sense, the passage from the quantale $\penq$ to a C*-algebra $\pena$, as that which is considered in this context by Connes~\cite{Connes}, represented quantalically by its spectrum $\Q \pena$, may be considered to correspond to the introduction of superposition. This view is further supported by the fact that each irreducible representation of $A$ has as Hilbert basis an equivalence class of Penrose sequences, suggesting that we may see each relational point of $\penq$ as being, in an appropriate sense, the Hilbert basis of a quantum point. However, further consideration of these ideas depends on a more careful examination of how $\penq$ is related to the C*-algebra $\pena$.

Most importantly, it is evident from the discussion throughout the paper that consideration of the quantale of Penrose tilings allows the introduction of a noncommutative space of which the relational points are exactly the equivalence classes of Penrose tilings of the plane. In the event that it may be shown that these are the only points of the quantale (and as has been remarked, none others are known), the quantale $\penq$ may be considered to represent the noncommutative geometric content of Penrose tilings in a particularly straightforward manner. It is to establishing this that we hope to return in a later paper.

\end{document}